\documentclass[a4paper,12pt]{article} 
\usepackage{bm,amsmath,amssymb}
\def\cat#1{{\mathfrak{#1}}}

\usepackage[latin1]{inputenc}
\usepackage[T1]{fontenc}

\newtheorem{making}{Making the Point}
\newtheorem{thesis}{Thesis}
\newtheorem{conjecture}{Conjecture}
\newtheorem{axiom}{Axiom}
\newtheorem{criterion}{Criterion}
\newtheorem{problem}{Problem}

\title{Halfway Up To the Mathematical Infinity I: \qquad
On the Ontological \& Epistemic Sustainability 
of Georg Cantor's Transfinite Design}

\author{Edward G. ~Belaga}

\begin{document}

\maketitle

\begin{abstract}

Georg Cantor was the genuine discoverer of the Mathematical Infinity, and whatever he claimed, suggested, or even surmised should be taken seriously -- albeit not necessary at its face value. Because alongside his exquisite in beauty ordinal construction and his fundamental powerset description of the continuum, Cantor has also left to us his obsessive presumption that the universe of sets should be subjected to laws similar to those governing the set of natural numbers, including the universal principles of cardinal comparability and well-ordering -- and implying an ordinal re-creation of the continuum. During the last hundred years, the mainstream set-theoretical research -- all insights and adjustments due to Kurt G\"odel's revolutionary insights and discoveries notwithstanding -- has compliantly centered its efforts on \emph{ad hoc} axiomatizations of Cantor's intuitive transfinite design. We demonstrate here that the ontological and epistemic \emph{sustainability} of th!
 is design has been irremediably compromised by the underlying peremptory, Reductionist mindset of the XIXth century's ideology of science. Our analysis and prompted by it synthesis lead to:
{\bf (i)} the extension of the well-known two-terms foundational opposition

\noindent$\bm{CN}:\ \{existence\ by\ axiomatic\ consistency \Rightarrow\ notational\ existence\},$

\noindent
to its novel, four-term axiomatic viability criterion $\bm{RSCN}$:

\noindent
$\ \ \{ontological\ relevancy\Rightarrow onto-epistemic\ sustainability\Leftrightarrow\bm{CN}\},
 $

\noindent
reducing $\bm{ZF}$ and its extentsions to the status of interactive programming languages manipulating \emph{ad hoc} contrived, pure notational infinite totalities,
{\bf (ii)} the new ontological insights into \lq\lq the nature'' of the continuum inspired by the quantum-mechanical entanglement argument, and {\bf (iii)} the interpretation of Cantor's class of all countable ordinals $\omega_1$ as an authentic, universal, ever emerging and never completed ordinal scale of the power and sophistication of iterative logical arguments.

\end{abstract}

\hfill\textit{I was beside the Master craftsman, }

\hfill\textit{delighting him day after day, }

\hfill\textit{ever at play in his presence, }

\hfill\textit{at play everywhere on his earth, }

\hfill\textit{delighted to be with the children of men.}

\vskip 2 mm

\hfill Proverbs 8:30-31

\tableofcontents

\section{Introduction}\label{intro}

\subsection{Cantor's Mission and His Peremptory Amalgams}\label{amalgam}

 \emph{Georg Ferdinand Ludwig Philipp Cantor}, 1835-1918, the genuine discoverer and the first forceful colonizer of the Mathematical Infinity's mountain range, was not just a scientist and mathematician: the powerful thinker of a strong spiritual bent, he was a missionary -- or, at least, so was he himself perceiving his scientific vocation, systematically pursuing it as the mission to shed a new intellectual and, hopefully, spiritual light into the mystery of the Infinite. 

It has never been easy to be a missionary, and as many other missionaries before and after him, Georg Cantor has fallen victim to the most sustained and cruel persecution -- in his case, academic. The rejection and the ridicule were so unforgiving -- especially, on the part of Leopold Kronecker, his elder and highly regarded colleague, the editor of the most prestigious German Mathematical magazine who \lq\lq\textit{considered Cantor a scientific charlatan, a renegade, a \lq corrupter of youth'}'' \cite{dauben1979} (p. 1) -- that Cantor's mental health has greatly suffered. He died in the mental institution unaware of the long overdue recognition his work finally received, in particular, from the London Mathematical Society.

There is no doubt that, alongside Cantor's health, the intellectual freedom, transparency, integrity, semantical relevancy and, ultimately, the ontological and epistemic sustainability of his scientific quest suffered as well. Most tragically, Cantor himself has been fully aware of this:

\lq\lq\textit{But despite the seeming thoroughness of it all, the special results of transfinite arithmetic did not compensate for the major flaw in Cantor's entire presentation {\rm[of his results in the major seminal paper of 1897 -- EB]}. The continuum hypothesis remained unresolved, as did the questions of whether every transfinite power was an aleph, whether the transfinite cardinals were all comparable, and whether every set could actually be well-ordered. Throughout the entire presentation, in fact, there was a nagging sense that something was not quite right.}'' (\cite{dauben1979}, p. 194)

Was it the psychological hardship aggravating the usual professional pressure to deliver the new ideas and results as fast and as strong as possible, or just the excitement to discover the new and infinite \lq\lq \emph{terra incognita}'', with its exquisite in beauty ordinal construction and its fundamental powerset description of the continuum? 

Whatever might be the reasons, the inescapable fact is that Georg Cantor has obsessively precipitated and stoically maintained the imposition of his original, overbearing global transfinite design, marked from our point of view by a grave ontological insufficiency, epistemic arbitrariness and driven by his Reductionist methodological presumptions, as follows: 

First, sets should be subjected to laws similar to those governing the set of natural numbers, including the universal principles of cardinal comparability and well-ordering  implying ultimately an ordinal re-creation of the continuum. This minor, suggestive, and relatively benign type of \emph{Cantor's reductionism}, with a lowercase \lq\lq\emph{r}'', has been exposed and analyzed on several occasions (cf., e.g., Michael Hallett's \cite{hallett1984}), albeit without any remedies proposed.

Second, the Universe of sets, or the Absolute, exists in the sense that is available to us as a whole for manipulation and, eventually, total formal control. 

This last, \emph{absolute type} of \emph{Cantor's reductionism} should not be confused with his or G\"odelian or, for that matter, Finslerian Platonism which is in principle compatible with an open, emerging universe of ideas about sets. Rather this militant Reductionism, forcefully imposing its pre-conceived, ontologically rigid and epistemicly arbitrary designs on the laws of the domain one is supposed \emph{to discover and not to invent}, was the product of the ideologically uniform, peremptory mindset of the XIXth century's philosophy of science (cf. Section~\ref{reductionism}). 

The result, as we see it, is most dispiriting: during the last hundred years, the mainstream set-theoretical research, all insights and adjustments due to Kurt G\"odel's revolutionary discoveries and insights notwithstanding, has compliantly and without great success centered its efforts on \emph{ad hoc} axiomatizations of Cantor's transfinite blueprint -- from Zermelo-Frenkel's Axiomatics to Axioms of Large Cardinals.

And yet, all such hindsight, backward reproaches notwithstanding, it is with the greatest delight and most sincere gratitude that we acknowledge here the importance for the present study of Georg Cantor's original optimistic phenomenological and ontological vision of the Mathematical Infinity. -- The vision intellectually most penetrating and epistemicly both most general from a mathematical practitioner's point of view and most generous from the point of view of a logician and a philosopher, has been one of the principal inspirations of the present study (Section~\ref{synthesis}). 

Our analysis and prompted by it synthesis will result in:

{\bf (i)} The extension of the well-known two-terms 
\lq\lq \emph{Platonic$\,\Rightarrow\,$formal\,}'' foundational opposition
$$\bm{CN}:\{existence\ by\ axiomatic\ consistency \Rightarrow\ notational\ existence\},$$
\noindent
to its novel, four-term axiomatic viability criterion (Sections ~\ref{platonic},~\ref{construal})
$$\bm{RSCN}:\ ontological\ relevancy\Rightarrow onto-epistemological\ sustainability\Leftrightarrow\bm{CN}\},
 $$
\noindent
reducing $\bm{ZF}$ and its extentsions to the status of interactive programming languages manipulating \emph{ad hoc} contrived, pure notational infinite totalities (Section~\ref{stumbling}).

{\bf (ii)} The new ontological insights into \lq\lq the nature'' of the continuum -- the insights inspired by the quantum-mechanical entanglement argument (Section~\ref{stumbling}).

{\bf (iii)} The interpretation of Cantor's class of all countable ordinals $\omega_1$ as an authentic, universal, ever emerging and never completed ordinal scale of the power and sophistication of iterative logical arguments (Section~\ref{redux}).

\subsection{Objectives and Results}\label{objectives}

Thus, one of the most important conclusions of our analysis will be the disavowal of Cantor's definition of the first uncountable ordinal $\omega_1$ -- the ordinal of the well-ordered set of all countable ordinals -- and the cardinal $\aleph_1$ associated with it  as a sheer notational gimmick void of any uncountable ontological substance  -- even if as such certainly helpful in notational designs of some \lq\lq very large'' countable ordinals \cite{sch\"utte1977}. 

To some degree, our argument here will parallel that of the non-existence of \emph{the ordinal of all ordinals}, famous because of its role in the foundational crisis in the beginning of the XXth century, and is inspired by the modern research on large countable ordinals \cite{pohlers1989}, with its systematic extension of our \lq\lq iterative upward mobility'' when searching for definitions of new, ever more \lq\lq fast growing'' countable ordinals -- leading to our interpretation of Cantor's collection, or class, $\omega_1$ of all countable ordinals as the ever emerging ordinal measure scale of the iterative power and sophistication of our formal logical endeavors (cf. Section~\ref{redux}). 

Deprived of its first uncountable ordinal, Cantor's ordinal theory of the continuum fails, with his Continuum Hypothesis retaining only its original, naif, non-ordinal interpretation (cf. Section~\ref{continuum}). 

Consequently, we are faced with the challenge to look for new, non-Cantorian ontological insights into \lq\lq the nature'' of the continuum -- the insights which could be eventually formalized by new set-theoretical axioms (cf. Section~\ref{locality}). 

Ours will be a still informal but far-reaching extension of the entanglement argument borrowed from quantum theory and the theory of quantum computation\cite{nielsen_chuang2000} (cf. Section~\ref{causation}). 

Armed with these novel transfinite insights and arguments, we abandon here for good both the framework of the classical set-theoretical reasoning firmly rooted, all protestations to the contrary notwithstanding, in the XIXth century's reductionist paradigm of mechanistic causality and the closely related to it -- in fact, implicitly underlying, even if postdating it -- classical theory of computation, formalized by Alfonso Church and Alan Turing (cf. Section~\ref{causation}). 

It is this theory which, in the absence of any formal alternative of a comparable importance, is deemed today by many to somehow stand for \emph{the theory of human thinking}, with Cantor's aforementioned \lq\lq\textit{fundamental law of thought}'' being a clear tributary to such a reductionist, outdated, but still beguiling philosophical appeal. Emboldened by our new insights into the Transfinite, we will conclude the present study by a discussion of an alternative approach to the essentially non-algorithmic, extemporaneous, creative human thinking (cf. Section~\ref{building}).

Some of our novel ontological and epistemic insights and suggestions are easily lending itself to appropriate formalizations in the today universally accepted Zermelo-Frenkel's set-theoretical framework suitably adjusted to our needs. 

This is the case of our understanding of the nature of Cantor's first uncountable \lq\lq ordinal'' $\omega_1$: its immediate formal definitions (cf. the notations in Section~\ref{redux})

$$\omega_1=\{\alpha\in On\mid\alpha\ is\ countable\}=\bigcup_{\alpha\ is\ a\ countable\ ordinal}\alpha
$$

\noindent
are either involving the General Comprehension Scheme with the condition \emph{countable} which, as we claim, is not even \emph{definite} (cf., e.g.,\cite{moschovakis1994}, p. 20) because involving the verification by counting to the infinite, or appealing to Zermelo's Fixed Point Theorem (\cite{moschovakis1994}, p. 73) coming at the heels of his well-ordering theorem which, in its turn, needs the full power of the Powerset Axiom and, for the sake of a formal justification of Cantor's construction of $\omega_1$, proves in our case too much to make this construction credible.

Other conceptual novelties of this study need new instruments of formalization which are at hand, too. Two particular tracks of such an eventual formalization are followed in our forthcoming papers \cite{belaga2009}, \cite{belaga_etal2009}.

We feel, however, that the full axiomatization of our blueprint of the Mathematical Infinity should be the subject of a separate presentation. This is why all formal aspects of our study are kept here to a minimum -- in consequence making eventually our paper accessible to a general mathematical reader.

\subsection{Bibliographical Note and Acknowledgments}\label{bio}

The historical remarks let drop above and below do not represent a history of either Cantor's discoveries, or its antecedents and sequels. Our choice is restricted to a few important for what follows conceptual breakthroughs whether, in our opinion, they are fully, or only partially, or not at all justified by the further developments \cite{belaga2008}. 

The interested reader could find mutually complementing accounts of the emergence and development of pre-Cantorian, Cantorian, and post-Cantorian Set Theories in the following papers and monographs listed according to the alphabetical order of their authors to whom the present author expresses his most sincere gratitude: 
Joseph Warren Dauben \cite{dauben1979},
Solomon Feferman \cite{feferman1987},
Michael Hallett \cite{hallett1984}, 
Thomas Jech \cite{jech2003},
Akihiro Kanamori \cite{kanamori2003}, 
Penelope Maddy \cite{maddy1988}, and
Gregory H. Moore \cite{moore1982}.

Author expresses his profound gratitude to these authors

\section{Epistemic Prequel to the Advent of Ordinals: Plain Continuum Standing Alone and High}\label{prequel}

\subsection{Prelude}\label{prelude}

It was originally the usually overlooked or intentionally ignored scandal of the undervaluation of the importance of the continuum, compared to that of Cantor's ordinal invention, as a priceless ontological source of the intuition into the Infinite that has struck a raw nerve with the author \cite{belaga1991}, \cite{belaga2001} and ultimately prompted the present conceptual revision. 

The epistemic and methodological roots of such a gross distortion, as well as the general historical and cultural predispositions to it are elucidated below, Sections~\ref{pricetags},~\ref{continuum}, as a prelude to our original approach to the foundational challenge posed by Mathematical Infinity -- the challenge whose meaning and importance were buried under the debris of diverse axiomatic and technical accommodations of Cantor's primeval transfinite design.

This is why, in the present chapter, we have chosen first to revise in a free, leisurely, and informal way some basic and vital concepts and ontological insights concerning the continuum, the concepts and insights coming both from Greeks and modern mathematics, but mostly forgotten, ignored, or hastily and uncritically absorbed by the Cantorian and post-Cantorian set theories -- if not irremediably crushed under the weight of their ever growing transfinite tower, with \lq\lq\textit{its top eventually reaching heaven}''. 

We accompany these revisions by informal and indicative for what follows hints of our principal insights, motivations, expectations -- briefly, of the driving intuitions behind the present study.

\subsection{Antiquity: an Infinity Flashback}\label{antiquity}

\hfill\textit{Heraclitus is the one who first} 

\hfill\textit{declared the nature of the infinite}

\hfill\textit{and first grasped nature as in itself}

\hfill\textit{infinite, that is, its essence as process.}

\vskip 2 mm

\hfill Georg Wilhelm Friedrich Hegel 

\hfill Vorlesungen \"uber die Geschichte der Philosophie

\vskip 6 mm

Before Georg Cantor has entered the scene of mathematical infinity, two types of infinite totalities -- the countable and the continuum, both re-discovered and radically re-mastered by Cantor -- were already well known \lq\lq experimentally'' to the mainstream Western mathematical and philosophical communities, starting with those of Ancient Greece. 

In fact, the Greeks have been the first to discover and \lq\lq colonize'' two basic infinite mathematical habitats, the set of \emph{natural numbers} denoted today $\mathbb N,$ or (Cantor's notations) $\omega$ or $\omega_0,$ and the \emph{real line continuum} $\mathbb R$ which have remained ever since both the most fundamental and widely used by modern practitioners of Mathematics, of Computer Science, and of Sciences at large, with their \lq\lq infinite sizes'', or \emph{powers}, according to Cantor, denoted respectively 
$\bigm\vert\omega\bigm\vert=\bigm\vert\omega_0\bigm\vert=\aleph_0$ and 
$\bigm\vert\mathbb R\bigm\vert=\cat{c}.$ 

These habitats existed for the Greeks as they exist for us, present-day mathematical yokels, independently of any formal or axiomatic justifications, but simply because they represent the indispensable, intuitively clear, and intellectually reliable mediums for beautiful mathematical theories rich in fruitful applications:

\lq\lq\textit{Among mathematicians, there is a widespread view that ongoing current mathematics on the whole is more reliable than any of the philosophically motivated programs that have been proposed to replace it, and that the only foundations that need be considered (if any at all) is organizational.}'' \cite{feferman2007}

Being well acquainted experimentally and fully comfortable theoretically with two basic numerical and geometric universes $\mathbb N$ and $\mathbb R$, the Greeks were fully conscious of the ontological differences between these two infinities, as is abundantly clear from the famous paradoxes advanced by Zeno of Elea (ca. 490 BC -- ca. 430 BC) (\cite{bochenski1970} p. 26; \cite{anglin_lambek1995} pp. 54-57). 

In particular, Zeno's paradox \lq\lq Achilles and the Tortoise'' clearly demonstrates the perceived by its author conceptual incompatibility between two different types of human experiences: the external, existential on the one hand and the inner, intellectual, and in particular mathematical, on the other hand -- the experiences nurturing two different concepts of infinity. 

The first type, associated today with the \emph{continuum} concept, emerges from our \emph{observations} of the external world viewed as an infinite and continuous one- or more-dimensional flow: the ever-changing skies, the waters of a river, the flight of an arrow, a running athlete, \emph{et cetera}. 

The second type, taking its origins in the human intellectual activity, is best encapsulated by all sorts of the \emph{counting experience} -- through observations of heart-beatings, walking as a step-by-step movement, building of towers, etc. -- and represents the only humanly available and never-ending \lq\lq accumulation of infinity'' by finite and discrete portions -- \emph{and so on}. 

Zeno, like Heraclitus before him, clearly doubted that the two infinities could be reconciled: one can run, but one cannot adequately understand or otherwise fully and exhaustively capture \emph{in a purely intellectual way}, and in particular formally, this phenomenon, because our reasoning proceeds by finite and discrete steps, by \emph{discrete markers} or \emph{buoys} occasionally dropped in the moving waters, whereas our and the river movement -- a mystery in itself -- is, as an ever-changing face of a river, continuous:

\lq\lq\textit{On those stepping into the same river, other and other waters flow} (Heraclitus of Ephesus, ca. 535 -- ca 475 BC) 

One should not be surprised to discover that Zeno's paradoxes remain of interest to modern philosophers as well. Thus, for instance, Henri Bergson (1859-1941) is thoroughly discussing and forcefully \lq\lq explaining'' them anew in his book \cite{bergson1959} (pp. 1259, 1376, 1377).

\begin{making}

\textbf{: Austere Ontological Realism and Epistemic Discretion of the Greeks.}
Heraclitus' and Zeno's visions are far from being Platonic. In their austere and minimalist epistemic way, they do not suppose the existence of an independently and objectively existing \lq\lq idea'' of \emph{the reality flow}, the idea which could be somehow and somewhen adequately and definitively captured by philosophical or mathematical formalizations. Instead, in the spirit of their austere ontological realism, they expect us to live in the presence of this flow -- the flow thought of as an ever ongoing, fully external to us, independent of us, and most rich phenomenologically process. The process to be patiently observed, experimented with, etc. -- with the expectation that we will be able with time to acquire new intellectual insights and intuitions somehow \lq\lq explaining'' this process, the intuitions susceptible to be eventually theorized and formalized. 

\end{making}

\subsection{From Formal to Platonic Existence, and Beyond.}\label{platonic}

It is clear that neither Georg Cantor, nor Kurt G\"odel, nor for that matter, Paul Finsler quoted below have shared such an extreme epistemic discretion and ontological \lq\lq price sensitivity'' -- choosing instead the Platonic vision of mathematical abstractions as objectively existing, directly accessible to us \lq\lq \emph{theoretical realities}'', without being initially inspired by prolonged and careful existential observations and experimental considerations, to be ultimately \lq\lq distilled'' from them as it were as pure abstract concepts and relationships:

\lq\lq\textit{Finsler's attitude towards mathematics was Platonistic in a very definite sense: he believed in the reality of pure concepts. Together they form the purely conceptual realm which encompasses all mathematical objects, structures and patterns. ... Mathematicians do not invent or construct their structures and propositions, they recognize, or discover, how these objects in the conceptual realm are interrelated with each other. It is clear that if there exists a conceptual realm, then it must be absolutely consistent; hence $\bm{existence}$ $\bm{implies}$ $\bm{consistence}$. ... The $\bm{Platonistic}$ $\bm{perspective}$ of mathematics can be expressed by the converse implication: $\bm{Consistency}$ $\bm{implies}$ $\bm{existence.}$ If a concept has been found to be consistent, it can be assumed to exist. This means that one can find properties and prove theorems about it.}'' (\cite{finsler1996}, p. 3)

Now comes the Formalist. Scarred by logical paradoxes which no Platonic force could exorcise, he is concerned only with the correct formal, notational existence and, respectively, formal deductibility of formal statements not leading ultimately to the contradiction $0=1.$ According to David Hilbert's axiomatic manifesto, axioms are not taken as self-evident truths: the elementary notions, such as point, line, plane, and others, could be substituted by tables, chairs, glasses of beer, etc. It is their defined relationships that only matter and are discussed. 

Or, as Abraham Robinson puts in the the infinity context which is here of main interest to us:

\lq\lq\textit{My position concerning the foundations of Mathematics is based on the following two main points or principles. $\bm{(i)}$ Infinite totalities do not exist in any sense of the word (i.e., either really or ideally). More precisely, any mention, or purported mention, of infinite totalities is, literally, meaningless . $\bm{(ii)}$ Nevertheless, we should continue the business of Mathematics \lq as usual', i. e., we should act as if infinite totalities really existed.}'' (\cite{robinson1965}, p. 230) 

Robinson's \lq\lq\emph{above the brawl}'' statement of \emph{the limited formalist liability} for any \lq\lq infinity risks'' undertaken by Mathematician both resembles and starkly contrasts with the infinity concept of Aristotle, which in its realistic, \lq\lq\emph{hands-on}'' approach surpasses in our opinion even the Platonic readiness to engage into the mathematical infinity combat. 

The following statement from Aristotle's Physics (Phys. III, 7, 207b27-34; the translation is borrowed from \cite{hintikka1966}, p. 201), anticipates the modern \lq\lq ultra-intuitionistic criticism'' of, say, Alexander Yessenin-Volpin \cite{yessenin1970} and explicitly raises the problem of the philosophical causality principles universally but implicitly underlying our mathematical queries (cf. below Section~\ref{causation}):

\lq\lq\textit{Our account does not rob the mathematicians of their study, by disproving the actual existence of the infinite ... In point of fact they do not need the infinite and do not use it. They postulate only that the finite straight line may be produced as far as they wish. It is possible to have divided in the same ration as the largest quantity another magnitude of any size you like. Hence, for the purposes of proof, it will make no difference to them to have such an infinite instead, while its existence will be in the sphere of real magnitudes.}'' 

As a matter of fact, it was Georg Cantor who has definitely demonstrated that \emph{mathematicians do need the infinite and they do use it} (cf., e.g., Section~\ref{redux}), notwithstanding the occasional expressions of doubt coming from our contemporaries -- from approvingly confident (\lq\lq\textit{Infinity in Mathematics: Is Cantor Necessary ?}'' \cite{feferman1987}) to panicky deserting (Abraham Robinson quoted above). 

The point is \emph{what do we need it for}, \emph{how do we use it}, and \emph{what might be the best workable and adequate philosophical and mathematical account of it}:

\begin{making}

: \textbf{Extending the Platonic-Notational Interplay into the Reality Context.} 

{\bf (1)} We express the interplay between Platonic and formal (notational) existence by the epistemic implication
$$\bm{CN}:\ \{existence\ by\ axiomatic\ consistency \Rightarrow\ notational\ existence\},
$$

{\bf (2)} It goes without saying that we do not deny, neither are we dismissive of, nor condescending to Platonic and formal mathematical deduction \lq\lq games'' which respect only one limitation and criterion of existence: the consistency of their ever extending axiomatic systems. Rather are we concerned here with another, and for that matter, crucial for us question: how far away and astray from the initial, fundamental aspirations and insights of our query into the infinity -- and ultimately, in what sort of a barren wilderness with no return -- are such games presently leading us \cite{belaga1998}, \cite{belaga2001}, Section~\ref{stumbling}. 

{\bf (3)} As Greeks before us, the majority of modern mathematicians could not be satisfied with the exclusivity of the epistemic implication $\bm{CN}$ pretending to dismiss as mathematically irrelevant or to ignore altogether the \emph{Ontological Relevancy}, the source and, through the \emph{Ontological and Epistemic Sustainability} link, the final judge of the \emph{Reasonableness}, i.e., \emph{Verifiable Explicative Strength} of object-oriented -- in this particular case, Continuum oriented -- mathematical theories.

{\bf (4)} Hence our four-term \lq\lq reality'' extension of the above epistemic implication $\bm{CN}$ (cf. also Section~\ref{construal}):
$$\bm{RSCN}: \{ontological\ relevancy\Rightarrow onto-epistemic\ sustainability\Leftrightarrow\bm{CN}\}.
$$

\end{making}

\subsection{Putting the Epistemic and Procedural Price Tags on the Cantorian Powerset Abstraction}\label{pricetags}

Georg Cantor has resolved (1874) one particular aspect of Zeno's \lq\lq Achilles and the Tortoise'' paradox, formally confirming Zeno's intuition, by introducing the method of \emph{one-to-one correspondence}, or \emph{bijection}, to both identify infinite totalities of the same \emph{power}, or \emph{equipotent} totalities, and to establish which of two totalities is \emph{bigger}. With this came his first, \lq\lq topological'' proof that the set of reals is uncountable.

Another one of Cantor's ingenious and powerful methods, that of \emph{diagonalization} (1891), has permitted him to prove even in a more direct and convincing way that \emph{the continuum} has the power strictly bigger than that of a \emph{countable} totality, $\cat{c}>\aleph_0$. 

Finally, Cantor has identified an operation (1891) acting on sequences of natural numbers, which \lq\lq re-creates'' the continuum from the totality of sub-totalities of the enumerable set $\mathbb N$ of natural numbers. In the finite case, the similar operation applied to a set with $n$ elements produces a set with $2^n$ elements -- hence the term \emph{the powerset construction,} $\cal{P},$ the one-to-one correspondents
$\cal{P}({\mathbb N})\cong\mathbb R,$
and the equality 
$2^{\aleph_0}=\cat{c}$.

Could all these original concepts of Cantor be understood by Greek philosophers and mathematicians? And would they agree that Cantor's powerset construction $\cal{P}({\mathbb N})\cong\mathbb R$ somehow invalidates Zeno's claim of the incompatibility of our two sources of intuition of the Mathematical Infinity, continuous and discrete? 

At any rate, it is clear that any relevant explanatory argument addressed to the Greeks should be adapted to their austere procedural realism: 

{\bf (1)} What do we mean by \lq\lq point'' of the continuum, the line, which is, for Heraclitus, Zeno and even Euclid, just a figure of an ever changing flow? How could one assign to this point a subset of natural numbers and what is the price of such an assignment?

{\bf (2)} \emph{Vice versa,} what do we mean by \lq\lq subset of natural numbers''? How could one assign to it a point on the continuum and what might be the price of \emph{this} assignment? 

{\bf (3)} And finally, what it means and how one might be sure that our procedures cover \emph{all} points and \emph{all} subsets ?

Now, even a brief analysis of Cantor's powerset construction applied to the set of natural numbers demonstrates that it does not invalidate but, quite to the contrary, vindicates Zeno's intuition and deepens its implications: 

{\bf (1)} The bijection $\cal{P}({\mathbb N})\cong\mathbb R$ is far from being symmetric and procedurally evenhanded. As a matter of fact, given a \emph{well-defined} subset of natural numbers, one is relatively free \lq\lq to mark a point'' on the continuum by constructing the corresponding real number with the help of, say, the continuum fraction device. 

{\bf (2)} To proceed, however, in the opposite direction, from a point \emph{well-defined} algebraically, analytically, or by another explicit, \lq\lq surfing over the continuum flow'' procedure -- distinct from the discrete in its nature powerset construction -- to a subset of natural numbers, one would need to execute an algorithm (such as the algorithm of Euclid calculating the corresponding continuous fraction from, say, a given system of equations) with a generally speaking infinite number of steps.

{\bf (3)} The \emph{irregular and unpredictable discontinuity character} of Cantor's powerset construction points in the same direction: the passage from the marker of a point on the continuum flow to the marker of another one is accidental, irregular, whereas the flow itself is moving along at every one of its points continuously and uniformly.

{\bf (4)} Moreover, we have no hope to ever reach \emph{all} points of this flow: 
\lq\lq\textit{not all propositions of the form
$$\alpha\ is\ a\ transcendental\ number
$$
are expressible, or definable in a language, since there are uncountably many transcendental numbers}'' (Finsler, 1926; cf. \cite{finsler1996}, p. 8.)

\begin{making}

\textbf{: Negligibility of the Impact of Our Occasional Stepping into the Continuum Flow.}

{\bf (1)} Alternatively to the \lq\lq\emph{stepping into a continuum flow}'', one can see Cantor's powerset device as sort of \lq\lq\emph{a scaffolding of the continuum skyscraper under perpetual construction}'', with each particular subset of natural numbers defining a transcendental number $\alpha$ being a discreet, sequential \lq\lq\emph{ladder of ascent}'', typically infinite, to $\alpha$ (cf. the epigraph to Chapter~\ref{ordinal_locality} below).

{\bf (2)} Our usage, above, of the qualifier \emph{well-defined} applied both to real numbers, or points of the continuum, and to (infinite) subsets of natural numbers is as intuitively transparent and informal as the one systematically employed by Cantor (cf., e.g., his \lq\lq law of thought'' statement in the above Section~\ref{amalgam}). 

{\bf (3)} It means that, without specifically restricting the interpretation of \emph{well-defined} to \emph{constructively defined}, according to one or another version of Constructivism, one can safely assume that the human activity, whatever might be its intensity, sophistication, and prolongation, might be able to actually \lq\lq probe individually into'' only countable many points of continuum and, respectively, subsets of natural numbers.

{\bf (4)} The extension of these apparently innocent observations about Cantor's powerset machinery of reaching out for points on the continuum to \emph{well-defined} subsets of the continuum will have substantial implications for Cantor's Continuum Hypothesis: cf. the next section. 

 {\bf (5)} The suggested limits of this machinery do not preclude the existence of other -- analytic, algebraic, number-theoretic, etc. -- methods of reaching out for the Uncountable. The interplay between such non-set-theoretic and set-theoretic methods is discussed in the next section.
 
 \end{making}

\subsection{Continuum Hypothesis: From Point Centered to Subset Centered Set Theory}\label{continuum}

Now comes Cantor's radical conceptual leap from the point-centered to the subset-centered approach to the theory of Continuum, and at this conceptual threshold -- leading, as David Hilbert claimed (\cite{hilbert1925}, p. 376), into \lq\lq\textit{the paradise Cantor created for us}'' from where \lq\lq\textit{no one shall be able to drive us}'' -- we have to abandon for good our Greek companions, albeit not their astute, austere epistemic methodology. 

To start with, Cantor has advanced (1878) the first, naif version of his \emph{Continuum Hypothesis}, which has appeared later as number one on the famous list of most important unsolved mathematical problems presented by David Hilbert during the Mathematical Congress at the University of Sorbonne, Paris, in 1900:
\begin{conjecture}

\textbf{$\bm{CHWE}$: Continuum Hypothesis in Its Weak Equipotency Form.} 
Given a well-defined sub-totality of the real line, $X \subset R$, it is either empty, or finite, or countable, or equipotent with the continuum.

\end{conjecture}

This \emph{weak equipotency} version of the Continuum Hypothesis, leaving to its prover to specify the subset of the continuum whose power she/he wants to compare with the four known powers, fits perfectly into Cantor's early mathematical experience (1872) of explicit construction of infinite sets of convergency points of trigonometric series. After introducing \emph{perfect sets}, Cantor was able to prove $\bm {CHWE}$ for all \emph{closed sets}. This was the beginning of Descriptive Set Theory \cite{kechris1995}

For better or worse, this has not been the end of the story. Being inspired by a much more radical vision of the discrete-continous opposition than Zeno, Cantor \lq\lq\textit{freed himself of all fetters and manipulated the set concept without any restriction}'' -- as Hermann Weyl puts it disapprovingly in \cite{weyl1949} (p. 50). In particular, Cantor's full equipotency version of the Continuum Hypothesis deals not with carefully specified sub-totalities, but either with \emph{any} sub-totality or with the collection of \emph{all} sub-totalities of the continuum:
\begin{conjecture}

\textbf{$\bm{CHFE}$: Continuum Hypothesis in Its Full Equipotency Form.} 
Any sub-totality $X$ of $\ \mathbb R$ is either finite, or countable, or equipotent with the continuum.

\end{conjecture}

This exercise assumes, implicitly or explicitly, the existence of the powerset of the continuum and thus extends the portability of the original Cantorian powerset construction beyond the countability case, where its legitimacy has been assured by the analytically defined bijection $\mathbb R\cong2^{\mathbb N},$ with the two \lq\lq mathematical habitats'' $\mathbb N$ and $\mathbb R$ well known \lq\lq experimentally'' and operationally long before the invention of the powerset gadget:

\lq\lq\textit{This does not prove the legitimacy of the} [universal] \textit{powerset principle. For the argument is not: we have a perfectly clear intuitive picture of the continuum, and the powerset principle enables us to capture this set-theoretically. Rather, the argument is: the powerset principle (or principles which imply it) was revealed in our attempts to make our intuitive picture of the continuum analytically clearer; in so far as these attempts are successful, then the powerset principle gains some confirmatory support.}'' (\cite{hallett1984}, p. 213.) 

\begin{making}

\textbf{: Cantor's Infinite Powerset Construction Is Fully Context-Dependable.}
The immediate and obvious extension of the above \lq\lq price-sensitivity'' analysis, Section~\ref{pricetags}, demonstrates the context-depen-dable character of Cantor's infinite powerset operation $\cal{P}(S), S$ being an infinite set: this operation is designed to act on, and only on the countable type of such sets, and it is meaningfully and ultimately formally sustainable only in this context. All other applications of this operation to infinite sets, as e.g., in the case of the \lq\lq set'' $\cal{P}({\mathbb R})$ of all real functions of one variable whose \lq\lq cardinality'' is denoted by $\cat{f}=2^\cat{c},$ are pure set-theoretical notations, certainly useful as such, but ontologically void of any verifiable reference to an \emph{objectively observable, meaningful reality context}.

\end{making}

Hence the following \emph{Thesis} -- the principle or conjecture which specifies only partially or informally, or leaves open altogether an important segment of its formal assumptions -- going back to Cantor's original insight, similar and, as we believe, ontologically and epistemically related, prior, and superior to Church-Turing's Thesis (cf. Section~\ref{causation}):
\begin{thesis}

\textbf{: On the Continuum Hypothesis.}

{\bf(1)} Whatever might be a specific extension of the qualifier \emph{well-defined}, any well-defined subset of the continuum is either empty, or finite, or countable, or equipotent with the continuum.

 {\bf (2)} The Continuum Hypothesis makes no sense outside the above \lq\lq well-defined'' version. 
 
 \end{thesis}
 
To give the reader a foretaste of possible conceptual realizations of the qualifier \emph{well-defined}, here is its simplest and popular \emph{ordinary mathematics} version: 
 
 \lq\lq\textit{We identify as} \emph{ordinary} \textit{or} \emph{non-set-theoretic} \textit{that body of mathematics which is prior to or independent of the introduction of abstract set-theoretical concepts. ... The distinction between set-theoretic and ordinary mathematics corresponds roughly to the distinction between \lq uncountable mathematics' and \lq countable mathematics.}'' (\cite{simpson1999}, p. 1) 
  
\subsection{Cantor's Unrestricted Uses of the Powerset Construction}\label{powerset}

Of course, Cantor himself was absolutely free from, in fact, blissfully unaware of the eventuality of any ontological and epistemic misgivings about unrestricted uses of his powerset device, the misgivings similar to those raised above. Quite to the contrary, the powerset operation was for Cantor not only universally applicable, but also unrestrictedly and eventually transfinitely iteratively auto-applicable:
\begin{axiom}

\textbf{: Powerset}

For each set $A$ there exists a set $B,$ called the \emph{powerset} of $A, \cal{P}(A),$ whose members, or elements, are the subsets of $A.$

\end{axiom}

Now, Cantor's diagonal argument used to prove that the cardinality of the continuum is strictly bigger than that of the set of natural numbers, Section~\ref{pricetags},
${\cat c}=\mit{card}(\mathbb R)=\mit{card}({\cal P}(\mathbb N))>\mit{card}(\mathbb N)=\omega,
$
is readily available for the proof of the general inequality $\mit{card}({\cal{P}}(A))>\mit{card}(A),$ which gives rise to Cantor's original transfinitely growing hierarchy of set:

$$\mathbb N, \mathbb R={\cal P}(\mathbb N), {\cal P}(\mathbb R),{\cal P}({\cal P}(\mathbb R)),\ldots
$$

\noindent
with their respective cardinalities

$$\aleph_0=\beth_0, \beth_1,\beth_2,\ldots
$$

It is at this \emph{first critical juncture} of the realization of his transfinite design, to assure that every set of the Universe of sets fits at some place in the emerging, powerset device driven hierarchy (with some less dramatic adjustments presently delivered by the Axioms of \emph{emptyset}, \emph{pairset}, \emph{extensionality}, \emph{union}, etc. \cite{moschovakis1994}) that Cantor needed his \emph{Universal Cardinal Comparability Principle}.

\subsection{Cantor's Universal Cardinal Comparability Principle and Two Ontologically Distinct Sources of Sets}\label{comparability}

>From the point of view of a naif set theory, the universal principle of cardinal comparability appeared to Cantor both immediate, natural, and ultimately more important than that of the his well-ordering \lq\lq law of thought'' to which we turn our attention in the next chapter. And Cantor's discovery (1882) of the comparative antisymmetry property
$A<_cB$ and $B<_cA\Longrightarrow A\sim B,$ 
known today as Schr\"oder-Bernstein theorem and being a part of what might be called \emph{ordinary mathematics} (cf. the above quote from \cite{simpson1999}), looked like just the first, and for that matter promising step in the right direction. 

However, as things stand today, the set-theoretical enforceability of the universal comparability principle demands the same axiomatic strength as, and is formally equivalent to Cantor's universal well-ordering principle (cf., e.g., \cite{moschovakis1994}, pp. 120-121). 

In other words, this enforceability came with a price tag which, more than hundred years after Ernst Zermelo's proof of this principle, is still hounding the research into the Infinite. 

It is instructive to trace the conceptual and formal breakthroughs leading to this state of set-theoretic affairs \cite{dauben1979}, \cite{moore1982}, \cite{hallett1984}. What is still missing, we believe, in such analyses is the understanding of the transformation of Cantor's original theory of the Infinite based on two, initially independent \emph{ontological sources of set tokens}, continuum and ordinal, into his global transfinite design which was steadily evolving ever since from favoring the foundational value of ordinals over that of the continuum to unreservedly -- and, as we claim, undeservedly -- absorbing the latter by the former.

Which brings us straight out to Cantor's ordinal invention.

\section{Ordinals: Their Awe-Inspiring Beauty, Their Necessary Uses, Their Peremptory Abuses}\label{awe_inspiring}

\lq\lq\textit{A team of Hollywood techno-wizards set out to} \lq bring 'em back alive'. \textit{So they took a little artistic license to make them half again as large. Anyway, what did books know? Then a surprising thing happened. In Utah, paleontologists found bones of a real raptor, and it was the size of the movie's beast.} \lq We were cutting edge'\textit{, says the film's chief modelmaker with a pathfinder's pride.}\lq After we created it, they discovered it'.'' (Andrea Dorfman, Behind the Magic of Jurassic Park, Times International 1993, n. 17, pp. 53-54.)

\subsection{The Beauty and Efficiency of Ordinal Construction Redux on the Outside the Cantorian Set Theory}\label{redux}

Independently of the listed above, Section~\ref{pricetags}, classical in their transparency and beauty clarifications of the fundamental relationships between two known to Greeks infinite mathematical habitats, $\mathbb N$ and $\mathbb R$, Georg Cantor had an extraordinary in its originality and beauty idea to extend the counting, or ordinal attributing, beyond natural numbers -- into the invented by him transfinite ordinal realm.

Became acquainted with Cantor's transfinite numbers, David Hilbert, for once, did not mince the words to praise it:

\lq\lq\textit{This appears to me to be the most admirable flower of the mathematical
intellect and in general one of the highest achievements of pure rational human
activity.}'' (\cite{hilbert1925}, p. 373)

Typically of David Hilbert, it was his enthusiasm, and not that of the original designer, Cantor (cf. the next section), that proved authentic, prophetic -- and totally selfless: 

{\bf (1)} Today, nobody is surprised that a research paper on, say, \emph{termination proof techniques for Term Rewriting Systems}, or \emph{TRS}, which play an important role in automated deduction and abstract data type specifications, starts as follows: 

\lq\lq\textit{Cantor invented the ordinal numbers}

\begin{equation}\label{ordinal1_dershowitz}
\begin{array}{lll}
0, 1, 2, 3, \ldots , n, n+1, \ldots \omega, \omega+1, \ldots \\
\omega2, \ldots \omega n, \ldots \omega^2, \ldots \omega^n, \ldots \omega^\omega, \ldots \omega\uparrow n, \ldots\\
 \varepsilon_0, \ldots \varepsilon_0^{\varepsilon _0} 
, \ldots \varepsilon_1, \ldots \varepsilon_{\varepsilon_0}, \ldots , and\ so\ on.\\
\end{array}
\end{equation}

 \noindent
 \textit{Each ordinal is larger than all preceding ones, and is typset of them all:} 
 
\begin{equation}\label{ordinal2_dershowitz}
\begin{array}{lllllllll}
\omega= the\ set\ of\ all\ natural\ numbers; \\
\omega2 = \omega\cup \{\omega+n\mid n\in \omega\}; \\
\omega n = \cup_{i < n}\omega i ; \\
\omega^2 = \cup_{n\in\omega}\omega n ; \\
\omega^\omega=\omega\uparrow2 = \cup_{n\in\omega}\omega^n ; \\
\omega\uparrow n = \cup_{i < n}\omega\uparrow i ; \\
\varepsilon_0 = \omega^{\varepsilon_0} = \cup_{n\in\omega}\omega\uparrow n ; \\
\varepsilon_0^{\varepsilon_0}= \omega^{{\varepsilon_0}^2}; \\
\varepsilon_1 = \cup_{n\in\omega}\varepsilon_0\uparrow n . \\
\end{array}
\end{equation}

\noindent
 \textit{The notation $\alpha\uparrow n$ represents a tower of $n\ \alpha s.$}'' (\cite {dershowitz1993}, p. 243)
 
Notice that in common parlance \lq\lq\textit{Cantor invented}'' or \lq\lq\textit{created}'', not \lq\lq\textit{discovered'}', the ordinal numbers.
 
 After the above most succinct and transparent introduction to ordinals, the author demonstrates how the \emph{ordinal descent} can be used to prove termination for particular classes of \emph{TRS}s, with the general \emph{TRS} termination problem (a special case of the halting problem for Turing machines) being of course undecidable. 

As to the ordinal descent, it is an important special case of descent along partially ordered sets, in particular, along trees. One of Cantor's most fruitful ideas has been the notion of a \emph{well-ordering}, \emph{WO}, i. e., of a linearly ordered set fulfilling the condition of \emph{finite descent, FD}, i. e., of termination after a finite number of steps of any descending subsequence (ordinals are, of course, special \emph{WO}s). The principal merit of the \emph{FD} condition is the extendibility of the mechanism of Mathematical Induction beyond natural numbers to any \emph{WO} and, in particular, to any ordinal. 

{\bf (2)} The first discovery of the \textit{necessary uses} (according to Harvey Friedman's favorite phraseology and the title of his seminal paper \cite{friedman1986}) of transfinite numbers in what has been called above \textit{ordinary mathematics} was made much earlier -- just ten years after Hilbert's pronouncement. 

Namely, Gerhard Gentzen \cite{gentzen1936} has proved that the \emph{validity of the law of mathematical induction} along Cantor's ordinal segment 
$0, 1, 2, 3, .\ldots \omega, \omega+1, \ldots \varepsilon_0$ , 
is equivalent to the statement of the \emph{consistency of Peano arithmetic}. It is worth to mention that, among other things, Gentzen's result has brought with it a dignified, even if only partial and, for that matter, countable rehabilitation of Hilbert's program \cite{rathjen2006}. 

Gentzen's approach gave birth to \emph{Ordinal Analysis} whose objective is to provide \emph{ordinal certificates} for the \emph{consistency} of more or less complete and constructive fragments of artihmetic and analysis: see the paragraph {\bf (5)} below.

{\bf (3)} About ten years after Gentzen, it was Reuben L. Goodstein who has discovered an elementary, natural, and yet number-theoretically meaningful and aesthetically appealing arithmetical proposition of Peano arithmetic \cite{goodstein1944}, readily understandable to high school students but not provable in Peano arithmetic -- because any proof of this proposition requires \textit{the necessary use} of a transfinite induction up to $\varepsilon_0$. 

Goodstein's has been also the first meaningful and elementary mathematical statement illustrating in Peano's elementary axiomatic framework the incompleteness theorem of Kurt G\"odel. Almost forty years later, Laurie Kirby and Jeff Paris, the re-discoverers of the misunderstood and forgotten result of Goodstein, have added to the emerging list of such elementary examples their \lq\lq Hercules against Hydra battle'' \cite{kirby_paris1982}.

{\bf (4)} Then, in 1949, Alan Turing has given a remarkable general interpretation 
of explicitly defined countable ordinals as succinct symbolic notations for algorithmic 
structures with multiple loops \cite{turing1949} -- the interpretation which inspired a series of remarkable results on the program verification illustrated, in particular, by the quoted above paper \cite{dershowitz1993} .

{\bf (5)} There is no other field of mathematics and mathematical logic where the suggestive, interpretive, and creative power of Cantor's ordinal construction is more pronounced and more effective than \emph{Ordinal Analysis} which builds on and integrates all aforementioned breakthroughs:

\lq\lq\textit{The central theme of ordinal analysis is the classification of theories by means of transfinite ordinals that measure their \lq consistency strength'
and \lq computational power'. The so-called} proof-theoretic ordinal \textit{of a theory also serves to characterize its provably recursive functions and can yield both conservation and combinatorial independence results.}'' (\cite{rathjen2006}, p. 45)

\begin{making}

\textbf{: Countable Ordinals as Formal Devices for Measuring the Onto-Epistemic Limits of Particular Iterative Methods.}

{\bf (1)} By their very meaning and formal definition, proof-theoretic ordinals should be countable \cite{pohlers1989}, \cite{sch\"utte1977}, which means that the study of more and more \lq\lq strong'' axiomatic theories goes hand in hand with the never ending upgrading of the hierarchy of \lq\lq large'' countable ordinals: 

$\bullet$ from $\omega$ to $\omega^2,$ the proof-theoretic ordinal of $RFA,$ rudimentary function arithmetic;

$\bullet$ $\omega^3$ and $EFA,$ elementary function arithmetic;

$\bullet$ $\omega^\omega$ and $PRA,$ primitive recursive arithmetic;

$\bullet$  $\epsilon_0$ and Peano arithmetic;
 
$\bullet$  $\Gamma_0,$ Feferman-Sch\"utte \emph{predicativity} ordinal, the proof-theoretical ordinal of $ATR_0,$ arithmetic transfinite recursion;

$\bullet\quad\ldots$ 

$\bullet$  and so on -- beyond predicativity, recursive notations, and eventually -- beyond all known today explicit notational and combinatorial descriptions, toward new, today unknown and unimaginable iterative insights and devices.

{\bf (2)} In particular, Ordinal Analysis provides us with \emph{the independent evidence of existence} of corresponding countable ordinals. In other words, \lq\lq\textit{after Georg Cantor has created them, Gerhard Gentzen, Reuben L. Goodstein and their followers have discovered them}''.

\end{making}

Hence the following \emph{Thesis} which is aspiring to capture the emerging in the above cases of Ordinal Analysis general pattern and to extrapolate it to future breakthroughs in our understanding of, paraphrasing \cite{rathjen2006} (p. 45)  the \lq\lq\textit{consistency strength, computational power, and combinatorial sophistication}'':   
\begin{thesis}

\textbf{: The Onto-Epistemic Interpretation of Cantor's Class of Countable Ordinals.}

{\bf (1)} The collection $\omega_1$ of countable ordinals is the authentic, universal, ever emerging and never completed formal ordinal measure scale of the power and sophistication of iterative logical arguments and methods. 

{\bf (2)} In particular, $\omega_1$ is \emph{a proper class}. 

{\bf (3)} The statement of such conceptual generality and formal vagueness cannot be either \emph{proved}, or \emph{falsified} otherwise than by a discovery of \emph{necessary use(s) of $\omega_1$ in ordinary mathematics providing an independent evidence of its existence}. 

\end{thesis}

{\bf Discussion. (1)} If one assumes the truth of the first part of the above thesis, the \emph{ontological justification} of its second part is readily and naturally available:

{\bf (2)} The assumption that $\omega_1$ is a set and, thus, ordinal implies roughly that one can create a new iterative method by uniting all iterative methods which are \emph{already in place or will be ever invented} -- which is meaningless. From the ontological point of view, the confusion comes from a substitution of implications of \emph{the real temporal eternity}, which could not be amplified, by those of \emph{mathematical infinity}, which permits different procedures of ever increasing \emph{iterative complexity}, to overstep all already discerned and formalized limits.

{\bf (3)} At present, there exists no convincing and sustainable results of necessary uses of the Uncountable. To give just one example: the method which relates the structure of elementary embeddings associated with large cardinals to that of self-distributive algebras and braid groups has been later supplanted by direct, not involving the Uncountable proofs \cite{buss_kechris2001}. This example illustrates both the methodological and proof-theoretical richness of the existing theories of the Uncountable and the eventual countable interpretability of the uses of such theories in ordinal mathematics.

We are elaborating the foundational implications of these informal insights below, Section~\ref{infinityabove}.

\subsection{Cantor's Transfinite Design Takes Shape -- at a Price: a Clear-cut Case of the Abuse of the Ordinal Device}\label{shape}

The above thesis, if true, disqualifies both Georg Cantor's definition of the first uncountable ordinal, Section~\ref{redux}, and the known axiomatic foundations of set theory, such as $\bm{ZFC}$ (Zermelo-Fraenkel with the Axiom of Choice, $\bm{AC}$), which count among their provable statements the well-ordering theorem and, thus, the ordinal characterization of $\omega_1$. The implications of this \emph{new foundational crisis} will be discussed in Chapter~\ref{infinityabove}.

But this will happen later ... For now, we still have to patiently follow the founder of set theory.

Similarly to what happened with the discovery of the powerset operation, Section~\ref{powerset}, it is at this \emph{second critical juncture} of the realization of his transfinite design, that Cantor -- not giving a damn about the only ordinal numbers he has actually designed, the \emph{countable} ones \cite{dershowitz1993}, and blissfully unaware of their future astonishing and brilliant \textit{necessary uses} \cite{friedman1986} -- has been exclusively interested in their eventual transposal into the domain of the uncountable. 

Ultimately, Cantor succeeded to do this by postulating that the \emph{well-ordered collection of all countable ordinals} is a set and, consequently, the \emph{least uncountable ordinal} $\omega_1$ of the cardinality $\aleph_1.$

With this methodologically peremptory, epistemicly perfunctory, ontologically unwarranted break-in into the uncountable, Cantor has finally fully consolidated his firm hold on the well-ordered design of the Universe of sets -- leaving to the future generations of mathematicians to gasp air for its viable axiomatic foundations and, among other things, to resolve the pure cardinal interpretation of the Continuum Hypothesis which suddenly became \emph{its only official version} (regretfully, even according to Kurt G\"odel's insightful presentation \cite{g\"odel1947}):
\begin{conjecture}

\textbf{: $\bm{CHCV}$: Continuum Hypothesis, the Cardinal Version.}

Prove that $\aleph_1=\cat{c}=2^{\aleph_0},$ or find another \emph{aleph} of the cardinal hierarchy equal to $\cat{c}.$

\end{conjecture}

The striking, almost number-theoretic in its simplicity ordinal-cardinal interpretation of the original phenomenological versions of the Continuum Hypothesis, Section~\ref{continuum}, became so attractive that virtually all speaking about the Continuum Hypothesis, such an authority in fundamental matters as $\bm {CH}$, Kurt G\"odel including \cite{g\"odel1947}, are formulating it as a conjecture about cardinals, $\bm{CHCV}$, without ever mentioning its more tangible origins.

Today, after more than hundred years of arduous efforts by thousands of researchers, some share a remote hope that Cantor might get it wrong and that, in fact,
$\aleph_1<\cat{c}=\aleph_2={\mit card}(\omega_2),$
with $\omega_2$ being the collection of all ordinals whose cardinality does not exceed $\aleph_1,$
or in other words, the ordinals which are either finite or countable or equipotent with the first uncountable ordinal $\omega_1 \ldots$
	
Others, and they are in the majority, either do not care anymore or disown, on foundational grounds of all possible persuasions, the very legitimacy of $\bm {CH}$, as in the following statement of Solomon Feferman \cite{feferman2007} squarely aiming at the cardinal version of $\bm {CH}$:

\lq\lq\textit{I came to the conclusion some years ago that $\bm {CH}$ is an inherently vague problem. This was based partly on the results from the metatheory of set theory showing that $\bm {CH}$ is independent of all remotely plausible axioms extending $\bm{ZFC}$, including all large cardinal axioms that have been proposed so far. In fact it is consistent with all such axioms (if they are consistent at all) that the cardinal number of the continuum can be \lq anything it ought to be', i.e. anything which is not excluded by K\"onig's theorem. The other basis for my view is philosophical: I believe there is no independent platonic reality that gives determinate meaning to the language of set theory in general, and to the supposed totality of arbitrary subsets of the natural numbers in particular, and hence not to its cardinal number.}'' 

\section{Posterior Axiomatizations of Cantor's Transfinite Blueprint: Foundational Challenge}\label{blueprint}

In this and next chapters, we are attempting to reconcile our ontological limitations on both the uses of the powerset device, Section~\ref{powerset}, and of the set-theoretic status of Cantor's ordinal collection $\omega_1$, Section~\ref{redux}, with the existing and intensively evolving since already hundred years theoretical and formal foundational tradition. 

To this end, we revise in the present chapter both {\bf (1)} the fundamental explicit and unspoken onto-epistemological and formal set-theoretical assumptions, the bedrock of $\bm{ZFC}$ axiomatics and its \emph{Large Cardinals} extensions, and {\bf (2)} its almost universally ignored sustainability implications.

\subsection{Ptolemaic-like Deadlock of $\bm{ZF}-$based and Ever Extending Axiomatics of the Cantorian Set Theory}\label{ptolemaic}

In the spirit of Kurt G\"odel's legacy of pure (meta-)mathematical and detached of any ontological considerations interpretations of formal axiomatizations of set theory, one cannot but conclude that Cantor's intuitive perception of the unlimited ascension along some \emph{platonically pre-existing ladder of large ordinals and cardinals} could be formally realized today only at the price of the invention of never ending extensions of the $\bm{ZF}$ axiomatic which are themselves belonging to some \emph{platonically pre-existing ladder of set-theoretical principles} ... 

And so on ? Indeed:

\lq\lq\textit{As our edifice grew, we saw how one by one the large cardinals fell into place in a linear hierarchy. This is especially remarkable in view of the ostensibly
disparate ideas that motivate their formulation. As remarked by H. Friedman, this hierarchical aspect of the theory of large cardinals is somewhat a mystery ... In other words, is there a hierarchy of set-theoretical principles in another galaxy above 
$\bm{ZFC}$, disjoint and incomparable to our large cardinals ?}'' 
(\cite{kanamory_magidor1978}, p. 104)

This paradigm of never ending axiomatic adjustments, at any given axiomatic juncture never sufficient to settle, say, Cantor's Continuum Hypothesis, reminds us of the Ptolemaic millenarian deadlock of unending and ever more complicated adjustments of his system of epicyclic planetary orbits to new, more precise astronomical observations of planets supposedly moving around the Earth -- \emph{the theological} (according initially to the Pythagorean tradition and later on to the Biblical one), \emph{ontological} (according initially to Aristotelean and then to scholastic philosophy), \emph{physical} and, in the line with all such religious, cultural, and scientific intellectual certainties, \emph{the mathematical} center of the observable Universe.

\subsection{Ontological and Epistemic Construals: Consistency, Relevancy, and Truth}\label{construal}

The eventual foundational repercussions of Feferman's \lq\lq\emph{above the brawl}'' pessimism (cf. the above quote from \cite{feferman2007}, the last section of the previous chapter) concerning the legitimacy of the Continuum Hypothesis, $\bm {CH}$, pale into insignificance compared to our actual disavowal of \emph{the onto-epistemological sustainability} of the whole Zermelo-Fraenkel axiomatics, $\bm{ZFC}$, and its \emph{Axioms of Large Cardinals} extensions -- from the applicability of the powerset device beyond its original countable limits, Sections~\ref{continuum}, ~\ref{powerset} (cf. Feferman's final remark), to the legitimacy of the ordinal characterization of the collection $\omega_1$ of countable ordinals, Sections~\ref{redux}, ~\ref{shape}.

Our \lq\lq\emph{hands-on}'' approach to the foundational challenge presented by Cantor's transfinite blueprint has everything to do with the four-term, \emph{onto-epistemological relevancy} scheme discussed in Section~\ref{platonic}: 
$$\bm{RSCN}:\ \{ontological\ relevancy\Rightarrow onto-epistemic\ sustainability\Leftrightarrow\bm{CN}\}.
 $$

Before putting our back into the foundational business, however, let us make it clear that our approach is not inspired by, or restricted to a particular philosophical or mathematical school of thought, be it Constructivism (even if we are quoting with an indubitable sympathy Hermann Weyl's constructivist puns), or Platonism (even if we fully sympathize with Kurt G\"odel's Platonic dictums), or something reasonably else. 

Our only criterion of \emph{the truth} of axiomatics whose consistency is either formally demonstrated or contextually implied is \emph{their explanatory power and success}, or \emph{onto-epistemological relevancy and sustainability}, as expressed by the above four-term $\bm{RSCN}$ scheme. 

Thus, in no way doubting the scientific and mathematical value of \emph{any} of the flourishing schools of modern set theory, we assert that any \emph{criterion of axioms plausibility} definitely and exclusively equalling \emph{consistency} to \emph{truth} is lacking substance and strength to secure \emph{the onto-epistemological relevancy and sustainability} of its eventual mathematical implications. 

The below \lq\lq\emph{consistency-equals-truth credo}'' of W. Hugh Woodin recognizes this by subjecting the truth of an axiom to, in this particular case, its ultimate \lq\lq\textit{number theoretic consequences}'' (\cite{woodin2003}, p. 31):

\lq\lq\textit{For me, granting the truth of the axioms for Set Theory, the only conceivable argument against the truth of this axiom} [Projective Determinacy Axiom]\textit{, would be its inconsistency. I also claim that, at present, the only credible basis for the belief that the axiom is consistent is the belief that the axiom is true. This state of affairs could change as the number theoretic consequences of the axiom become more fully understood.}'' 

\subsection{Constative $\bm{vs}$ Performative Axiomatic Paradigms: Performative-Iterative Wishful Axiomatic Thinking}\label{wishful}

Our main foundational argument, developed below, Section~\ref{wishful} concerns the nature of axioms of Set Theory, from $\bm{ZF}$ to $\bm{ZFC}$ to \emph{Axioms of Determinacy} to \emph{Axioms of Large Cardinals}, which are typically \emph{inductive principles and procedures of emergence} of set-theoretical entities acting on swaths of transfinite ordinals and cardinals \cite{jech2003}. 

This is why they postulate the characteristics of such entities at least as much as -- and probably, transfinitly less than -- they postulate our computational and deductive proof-theoretical credentials -- \emph{finite, countable, and uncountable} \cite{pohlers1989}, \cite{rathjen2006}.

This is also why are misplaced all analogies of set-theoretical axiomatics with, say, geometric axiomatics, as are misplaced the analogies between functions defined directly and by finite procedures -- analytically, algebraically, or number-theoretically -- and algorithmically, i.e., by \emph{a priori} arbitrary recursive functions and procedures.

Here is just one and, for that matter, well-argued example of the popular genre of set-theoretical apologetics:

\lq\lq\textit{For instance, the proof that the axiom of parallels does not follow from the other Euclid axioms did not close geometry, but made the emergence of non-Euclidian geometries possible, and opened the question of recognizing, among all possible geometries, the most relevant for describing the physical world. Likewise, G\"odel's and Cohen's results show that several universes are possible from $\bm{ZFC}$, and, therefore, they open the study of the various possible universes.}'' \cite{dehornoy2003}

As a matter of fact, the above analogy only sharpens up our case: 

{\bf (1)} The Axioms of Euclidean geometry do not \emph{generate} objects, and being strictly fact-finding, or \emph{constative}, not creative or \emph{performative}, utterances, as Set-Theoretical Axioms, they were carefully distilled by Greek mathematicians from every-day geometrical experiences, representing straightforward, almost banal in their simplicity and transparency abstractions of basic and immediately and universally observable objects and their relationships. 

{\bf (2)} As to Set-Theoretical Axioms, they are designed mostly \emph{ad hoc} to urgently and in many cases rather arbitrary fill in the logical and iterative lacunas and thus to justify \emph{after the fact} some risky iterative flights of our unruly mathematical fancy. In short, those are axioms \emph{one needs to believe} \cite{maddy1988} rather than be banally sure of.

{\bf (3)} As to \lq\lq\textit{the emergence of non-Euclidian geometries becoming possible}'' -- in the 19th century, two thousand years after the Euclidian axiomatics -- this historical fact should be in the first place attributed to the general cultural trend inspired and guided by the accumulation of the scientific experience related to \emph{the geometric relativity} and closely associated with the slow replacement of the Ptolemaic geocentric system (emerged in the first half of the 2nd century to be later dogmatically asserted in the Western Christendom) by Copernicus' heliocentric system made public in 1543. 

{\bf (4)} This does not mean that we doubt the portability of learning experiences gained as the result of the passage from the strictly Euclidian to non-Euclidian geometries. Quite to the contrary: we believe that there is an important lesson to learn from this historical experience -- any true revolution in mathematics comes as a fruit of an acquisition of radically novel ontological and epistemic, which means also universal cultural and intellectual, paradigms.

{\bf (5)} It is also true that, in contradistinction to geometry,  the origins of modern set-theory were unashamedly metaphysical and even \lq\lq theological'', as are today the inspirational impulses of new axiomatic initiatives. 

Cantor himself was quite unapologetic about his motives. Here he is, writing hundred years ago to Father Thomas Esser in Rome (\cite{meschkowski1964}, p. 94):	
	
\lq\lq\textit{The establishing of the principles of mathematics and the natural sciences is the responsibility of metaphysics. Hence metaphysics must look on them as her children and as her servants and helpers, whom she must not let out of her sight, but must watch over and control, as the queen bee in a hive sends into the garden thousands of industrious bees, to suck nectar from the flowers and then together under her supervision, to turn it into precious honey, and who must bring her, from the wide realm of the material and spiritual world, the building blocks to finish her palace.}'' 

Today, philosophical and metaphysical insights continue to shape -- both implicitly and explicitly and more than any other source of formal intuition -- the bulk of \lq\lq theological ventures'' of modern set theory:

\lq\lq\textit{The adaptation of strong axioms of infinity is thus a theological venture, involving basic questions of belief concerning what is true about the universe. However, one can alternatively construe work in the theory of large cardinals as formal mathematics, that is to say the investigation of those formal implications provable in first-order logic.}'' (\cite{kanamory_magidor1978}, p. 104)

\subsection{Cantorian Operational and Generative and Hilbertian Universal and Argumentative Reductionism}\label{reductionism}

\textbf{Alert to the Intended Meaning of the Reductionism Dismissal.} 

\emph{The focal point of the following critical remarks is the dangers of the reductionist dismissal of ontology in the favor of pre-programmed ideological imperatives -- one of the main topics of the philosophical interest of the present study. Our criticism does not concern the reductionist methodology by itself, which is an indispensable and powerful instrument inseparable from the basic scientific idea of breaking a difficult problem into much more amenable to formal treatment pieces.} 

All eventual extraneous criticisms notwithstanding, it was and remains \emph{the reductive enemy within} which is most harmful to the free exercise of a researcher's metaphysical imagination. Thus, Albert Einstein, the discoverer of physical quants, could not accept non-local implications of the Quantum Theory because of his insistence of the universality of local causality.

The \emph{absolute type} of \emph{Cantor's reductionism} should not be confused with his or G\"odelian or, for that matter, Finslerian Platonism which is in principle compatible with an open, emerging universe of ideas about sets. 

Rather this militant Reductionism, forcefully imposing its pre-conceived, ontologically rigid and epistemicly arbitrary designs on the laws of the domain one is supposed \emph{to discover and not to invent}, was the product of the ideologically uniform, peremptory mindset of the XIXth century's philosophy of science. 

Other examples abound. In fact, one can assert somewhat schematically that:

{\bf (i)} the diversity of the material Universe \emph{was reduced} by Simon-Pierre Laplace to determinist implications of the laws of Newtonian Mechanics, 

{\bf (ii)} while the destiny of civilizations, according to Karl Marx, is governed by economic laws and struggles of classes, 

{\bf (iii)} with the natural section mechanism of Charles Darwin being made responsible for the diversity, riches and beauties of the biological life, 

{\bf (iv)} whereas Sigmund Freud was proclaiming the sexual drive, \emph{libido} to be the single fundamental creative force behind the personal destiny of man.

In the set-theoretical context, it was this local generative reductionist liability that has most distorted, in the present author's opinion, the future development of Set Theory -- by insisting on perpetual inventing and piling up new transfinite generative mechanisms as the only remedy to persisting explanatory lacunas. 

Two of these mechanisms invented by Cantor himself stand out  as apparently immediately \lq\lq given'': the powerset principle and the accumulation of new alephs by transfinite induction coming from two different and, as we claim, Chapter~\ref{locality}, ontologically distinct and epistemicly unrelated sources of mathematical experience, the continuum and the discrete. 

To summarize our grievances, we contest here, first, the blind, automatic extrapolation and extension of these two procedures beyond their original ontological matrices, respectively: the powerset representation of the continuum and the universe of countable ordinals. Similar \lq\lq\textit{more is different}'' doubts concerning the reductionist ignoring the scale of the viability and the emergent nature of established concepts and laws have been raised by physicists \cite{anderson1972}, \cite{laughlin_pines2000}. Second, we contest the ontological soundness of the amalgamation of these two procedures into one \lq\lq linear ascending ladder'' (cf. Cantor's formula in the epigraph to Chapter~\ref{ordinal_locality} below) of $\aleph$'s and $\beth$'s, Section~\ref{powerset} ...

... And then, in the middle of the first foundational crisis, it was the turn of David Hilbert to enter the fray with his \emph{global argumentative reductionism} which was attempting to reduce all mathematical reasoning to a fully formalized mechanical proof-theoretic procedure. 

Until, that is, Kurt G\"odel's has demonstrated that the unexpected and inexplicable are always there to defy our \emph{currently available} fully formalized and assertive self-confidence \cite{zach2006}.

\subsection{Hilbertian, Post-Hilbertian, G\"odelian Programs and Their G\"odelian \& Post-G\"odelian Stumbling-Blocks}\label{stumbling}

David Hilbert's foundational expectations and their brutal refutal by Kurt G\"odel created the most disheartening paradigm of the formal scientific thinking: the stronger are the axiomatic foundations or, in other words, the expressive power of a formal theory, the more inevitably stumbles one upon very natural and formally admissible in this theory claims which can be neither proved nor refuted without some nontrivial extension of the said axiomatic foundations. And in most interesting cases, the extension in question has to be either the claim itself or -- apparently, rather arbitrary and with equal plausibility -- its logical negation ! \cite{belaga1998}, \cite{belaga2001}

The challenge to recover in this atmosphere of consistency relativism the viable foundations of set theory has been ultimately and, for that matter, most reluctantly accepted by Kurt G\"odel. 

Unfortunately, the liberator of the foundational metamathematical research from the illusions of the Hilbertian argumentative reductionism fell himself under the spell of the Cantorian ideological Reductionism casting itself as Platonism. Most regrettably, discussing the importance and perspectives of the Continuum Hypothesis, $\bm{CH},$ in his seminal article \cite{g\"odel1947}, Kurt G\"odel has wholeheartedly accepted Cantor's transfinite design and the cardinal interpretation of $\bm{CH}.$ This inexorably led him into advancing a program of the search for new \emph{performative} axioms of set theory -- the axioms that are strong enough to answer questions left undecided by the standard axioms $\bm{ZFC}.$ 
 
 G\"odel's reputation for perspicacity, prudence, and intellectual integrity has been apparently decisive in influencing the ensuing \emph{gold rush} into \lq\lq High Infinite'' \cite{kanamori2003}. 
 
 The excesses of this very abstract, almost ideological trend could be corrected by a thinker of a more realistic, computational cast, such as Alan Turing who discovered both an alternative, algorithmic interpretation of countable ordinal procedures \cite{turing1949} (cf.  Section~\ref{redux}), and the \emph{Halting Barrier},  the undecidability of the Halting Problem, an algorithmic projection of G\"odel's incompleteness results: there is no \lq\lq halting Turing machine'' capable of distinguishing between halting and non-halting programs \cite{turing1936}. In other words, given a program, the only way to discover whether it halts or not is to run it -- possibly, indefinitely.
 
 Today, sixty years after the publication of the first installment of G\"odel's program \cite{g\"odel1947}, here is an informal version of such a correction:
 
 \begin{thesis}
 
 \textbf{: Post-Turing Halting Barrier for Performative Set-theo-retical Axioms of Iterative Nature.}
 
There exists neither general metamathematical principle, nor logical criterion, nor verifiably terminating computational procedure to establish the objective and substantial \lq\lq truth'' of a performative set-theoretical axiom of iterative nature which postulates the existence of a transfinite object outside the already existing (say, $\bm{ZF}-$based) transfinite scale -- otherwise that is than \lq\lq to run'' the theory completed with the new axiom until it would be discovered some independent \lq\lq necessary uses'' of the object in question. 
 
 \end{thesis}

There is no other modern domain of formal studies where this post-modern paradigm is so pronounced as in Set Theory and its $\bm{ZF}$-based axiomatics. The trouble is hidden exactly where our $\bm{ZF}$ pride resides: in the powerful built-in iterative mechanisms of set generation. In other words, $\bm{ZF}$ has gained in its creative power on the expense of its descriptive power, becoming a sophisticated \emph{programming language}, which is successfully mimicking some aspects of the Mathematical Infinity but whose main thrust lies with the providing to advanced \lq\lq users'' sophisticated options of generation of, and manipulation with artificial transfinite totalities, similarly to computer graphic imagery of video games -- with {the Axiom of Determinacy} opening the advent of \emph{Infinite Games} \cite{kanamori2003}. 
	
Thanks to this interpretation, finds its proper place, in particular, the puzzling and disconcerting predominance in modern set theory of results on $\bm{ZF}$ consistency and independency: 

\lq\lq\textit{When modern set theory is applied to conventional mathematical problems, it has a disconcerting tendency to produce independence results rather than theorems in the usual sense. The resulting preoccupation with \lq consistency' rather than \lq truth' may be felt to give the subject an air of unreality.}'' (\cite{shelah1992}, p. 197) 

Finally, we claim that all such consistency results are just the instances of successful program verification. 

In other words, the totalities in question are, in fact, pure mathematical notations not related to any \lq\lq set-theoretical reality'' outside the tight structure of their definitions and relationships which might turn out to be suggestive of our permanently evolving iterative programming ability: 

\lq\lq\textit{Only the first few levels of the cumulative hierarchy bear any resemblance to external reality. The rest are a huge extrapolation based on a crude model of abstract thought processes. G\"odel himself comes close to admitting as much.}'' (\cite{simpson1988}, p. 362)

These \emph{a fortiori} set-theoretical observations and related logical heuristics suggest the following informal thesis which projects G\"odel's incompleteness theorems into the realm of the emerging art of \lq\lq mathematical novelization'': 

\begin{thesis}

\textbf{: On the Post-G\"odelian Incompleteness.}

Any conceptually sufficiently rich and logically/axiomatically sufficiently sophisticated mathematical system allows a huge, super-exponentially expanding \lq\lq\emph{mathematical sci-fi novelization}'' -- the creation of a multitude of \lq\lq\emph{mathematical sci-fi novels}'', i.e., \emph{fully consistent} mathematical theories with unlimitedly extending axiomatic bases -- \emph{\lq\lq forced themselves upon us}'', as it were, not only \emph{\lq\lq as being true'' (\cite{g\"odel1964}, p. 268)} but being also intellectually compelling and esthetically attractive -- and yet which do not have in their (more than) overwhelming majority, either at this juncture or whenever in future, any objectively verifiable mathematical and/or substantial extra-mathematical meaning outside the proper, self-absorbed scene of formal deductions inside the system in question.

\end{thesis}

\begin{criterion}

\textbf{Sustainability Criterion of Viability.}

{\bf (1)} This new mathematical reality suggests the following  \emph{Onto-Epistemic Sustainability Criterion of Viability}: an axiomatic should ultimately stick as close as possible to the onto-epistemologically sustainable notions and relationships of the underlying mathematical domains:
$$\bm{RSCN}:\ ontological\ relevancy\Rightarrow onto-epistemological\ sustainability\Leftrightarrow\bm{CN}\},
 $$

{\bf (2)} The qualifier \lq\lq ultimately'' to the above assertive \lq\lq stick to'' is here not only to intimate that our \emph{Criterion of Viability} is unenforceable and, thus, voluntary. It is also a warning justified by the current mathematical experience, and not only set-theoretical, that axiomatic efforts in a theory with a rich ontological base could produce attractive and apparently plausible conjectures and problems of monster logical and proof-theoretical complexity, with the daunting task to prove or disprove them amounting to the mobilization of intellectual efforts of all humanity for some centuries. 

{\bf (3)} The examples of the Ptolemaic astronomic system and of the famous unresolved mathematical problems of Greeks confirm such an eventuality -- even if on a conceptually much more modest scale. It means that the mathematically \lq\lq free-market, easy-going'' attitude toward the choice of, and the level of human and institutional investments into the mathematical domain should be somehow assisted by responsible and trusted \lq\lq oracles''. This said, our pundits could be also wrong, as we believe was wrong Kurt G\"odel directing us into the Large Cardinals morass.

\end{criterion}

\section{Nostalgic Interlude. From the Infinity above to the Infinity below: Continuum \& Ordinals}\label{infinityabove}

In his spiritual autobiography, the famous Russian writer and thinker Lev Nikolayevich Tolstoy transcripts his night dream, somewhen in 1879, of lying on his back in a bed and trying \lq\lq\textit{to consider how and on what was} [he] \textit{lying -- a question which had not till then occurred to}'' him \cite{tolstoy1882} (pp. 46-47):

\lq\lq\textit{And observing my bed, I saw I was lying on plaited string supports attached to its sides. ... I seemed to know that those supports were movable. ... I made a movement with my whole body to adjust myself, fully convinced that I could do it at once; but the movement caused the ... supports under me to slip and to become entangled, and I saw that matters were going quite wrong. ... I looked down and did not believe my eyes. ... I could not even make out whether I saw anything there below, in that bottomless abyss over which I was hanging and whither I was being drawn. My heart contracted, and I experienced horror. ... What am I to do ? What am I to do ? I ask myself, and look upwards. Above, there is also an infinite space. I look into the immensity of sky and try to forget about the immensity below, and I really do forget it. ... I know that I am hanging, but I look more and more into the infinite above me and feel that I am becoming calm. ... And then, as happen!
 s in dreams, I imagined the mechanism by means of which I was held; a very natural intelligible, and sure means, though on one awake that mechanism has no sense.}''

Compare this with Hermann Weyl's final remarks to his short and lucid review of the foundational efforts set in motion by the discovery of logical paradoxes, \textit{the challenge which had not till then occurred to mathematicians}:

\lq\lq\textit{From this history one thing should be clear: we are less certain than ever about the ultimate foundations of (logic and) mathematics. Like everybody and everything in the world today, we have our \lq\lq crisis''. We have had it for nearly fifty years. Outwardly it does not seem to hamper our daily work, and yet I for one confess that it has had a considerable practical influence on my mathematical life: it directed my interests to fields I considered relatively \lq\lq safe'', and it has been a constant drain on my enthusiasm and determination with which I pursued my research work. The experience is probably shared by other mathematicians who are not indifferent to what their scientific endeavours mean in the contexts of man's whole caring and knowing, suffering and creative existence in the world.}'' (\cite{weyl1946}, p. 13)

Today, hundred years after the emergence of the first foundational crisis in Mathematics, an ordinary mathematician surfs fearlessly on an opportune wave of new conceptional and formal breakthroughs, carrying out routinely all necessary to him operations, whether they are supposed to be \textit{carried out by God} or by men. 

And yet, his carelessness does not affect in the slightest the reality and importance of the abyss below his feet -- \emph{Ordinal transfinite}, Chapter~\ref{awe_inspiring} -- and the immensity of the infinity above his head -- the mathematical habitat called \emph{Continuum}, Chapter~\ref{prequel}, associated with, and deriving its permanently renewing onto-epistemological relevancy and suggestive power from our being active \lq\lq\textit{in the presence of the reality flow -- the flow thought of as an ever ongoing, fully external to us, independent of us, and most rich phenomenologically process}'' (Section~\ref{antiquity}).
\begin{thesis}

\textbf{: Two Sources of Mathematical Infinity.}

{\bf (1)} Since at least three thousand years, the Countable and the Continuum represented two sources of the mathematical intuition associated today with Mathematical Infinity. 

{\bf (2)} All set-theoretical axioms and infinity notations notwithstanding, the exclusivity of these sources has been confirmed by the mathematical developments of the last century.   

{\bf (3)} There exist no other ontological sources of Mathematical Infinity. This informal claim complements and generalizes our Continuum Hypothesis Thesis, Section~\ref{continuum}.

\end{thesis}

\section{Local Causation of  \lq\lq Man's Mathematics'' 
\emph{Versus} Non-Locality of Classical Mathematics}\label{locality}

\hfill\textit{Our point of view is to describe the mathematical\ }

\hfill\textit{operations that can be carried out by finite beings,}

\hfill\textit{man's mathematics for short.\qquad\qquad\qquad\qquad\quad}

\hfill\textit{In contrast, classical mathematics concerns itself }

\hfill\textit{with operations that can be carried out by God.\ \quad}

\vskip 2 mm

\hfill Errett Bishop \cite{bishop1985}, p. 9.

\vskip 6 mm

\subsection{The Continuum, Suslin's Problem, Local Causation, Quantum Non-locality, and Church-Turing Thesis}\label{causation}

And what about other, \emph{classical and \lq divine' Mathematics} (as Errett Bishop's has chosen to descried it)? How about Mathematics of a free surfing on the real line continuum, starting with the logically unimpeded Classical Analysis and passing by the famous and still open Suslin's conjecture?
\begin{problem}

\textbf{: Suslin's Conjecture.}

Let ${\mathbb K}$ be a linearly ordered set without the first or last element, connected in order topology, with no uncountable family of pairwise disjoint open intervals. Is ${\mathbb K}$ isomorphic to the real line ${\mathbb R}$?(\cite{kunen1980}, p. 66)

\end{problem}
\begin{thesis}

\textbf{: Non-Locality Characterization of the Continuum.}

{\bf (1)}The striking feature of  classical analytic machinery conceived to deal with, and perfectly adapted to the Continuum habitat, and which immediately distinguish it from methods and theories subject to the discrete, Ordinal Analysis related treatment, is its intrinsic, inextinguishable, fundamental, outside of the Continuum not existing and not obtainable \emph{non-locality}, in the sense this term is understood in \emph{Quantum Information Processing}, or \emph{QIP}, for short \cite{nielsen_chuang2000}.

{\bf (2)} Any phenomenologically and ontologically faithful or at least relevant axiomatization of the Continuum should include a \emph{Non-locality Postulate}, or \emph{Non-locality Axiomatic Scheme}, to formally account for the following property of the Continuum: 

All \lq points', or \lq elements' of the Continuum are non-locally, i.e., simultaneously and at any moment, accessible. This non-local accessibility extends to all well-defined \lq slices' (subsets in the set-theoretical terminology) of the Continuum.

\end{thesis}

{\bf Comments and Implications: From Aristotle to Suslin to Church-Turing Thesis. (1)} The quoted above claim of Aristotle about the mathematical sufficiency of \lq\lq\textit{the sphere of real magnitudes}'' explicitly denies among other things such a non-locality characterization of the Continuum. One can only admire the cleverness, unambiguity, and consistency of Aristotle, the father of formal logic, who was sticking on this occasion to his and his colleagues, from Zenon to Euclid, \emph{discrete, local causality interpretation of mathematics}. For the sake of completeness and facility of the understanding our argument, we reproduce here Aristotle's dictum (see for details Section~\ref{platonic}):

\lq\lq\textit{In point of fact they} [mathematicians] \textit{do not need the infinite and do not use it. They postulate only that the finite straight line may be produced as far as they wish. It is possible to have divided in the same ration as the largest quantity another magnitude of any size you like. Hence, for the purposes of proof, it will make no difference to them to have such an infinite instead, while its existence will be in the sphere of real magnitudes.}'' 

{\bf (2)} We believe that, the availability of a Non-Locality Postulate in a modified $\bm{ZF}$-like axiomatics -- our powerset and ordinal limitations of size being assured, Sections~\ref{powerset},~\ref{redux} -- Suslin's conjecture becomes an easily provable statement (we will be back to this claim in our forthcoming paper).

{\bf (3)} To make these observations amendable to an eventual formalization, Section~\ref{ordinal_locality} below, let us turn our attention to Robert Gandy's well-known real-life analysis \cite{gandy1980} of Church's and Turing's computability, and in particular, to Gandy's \textit{Principle of Local Causation}, as it is informally and succinctly summarized by David Israel (\cite{israel2002}, p. 197):

\lq\lq\textit{Causal effects must be locally propagated ... There is an upper bound (e.g., the velocity of light) on the speed of propagation of changes. ... There is no (unbounded) action at a distance; no simultaneous causation.}''

{\bf (4)} After adding to these principle three other, more material and routine \emph{Principles of Mechanism}, Gandy demonstrates that, first, functions computable on such abstract devices are 
\lq\lq\textit{simply the Turing-computable functions, thus adding a striking bit of evidence for the adequacy and stability of Turing's analysis}'' (\emph{ibid.}), and second, that each of these four conditions is necessary to avoid either absurdity, or -- especially, in the case of the Principle of Local Causation -- indiscriminate, over-reaching applicability and not restricted to recursive functions computability.

{\bf (5)} In a similar vein, Nachum Dershowitz and Yuri Gurevitch \cite{dershowitz_gurevich2008} have recently identified three interrelated \emph{Sequential Postulates} \cite{gurevich2007} which represent an algorithmic axiomatization of computability allowing for a proof of Church's and Turing's theses, with the first of these postulates representing an algorithmic equivalent of Gandy's Principle of Local Causation.

\subsection{Ordinal Constructions Abide by Principles of Local Causation}\label{ordinal_locality}
 
\hfill\textit{Following\ the\ \ finite there\ is a} transfinite, \textit{that is}

\hfill\textit{an unbounded ascending ladder of definite modes,}
 
\hfill\textit{which by their nature are not finite\ \ but\ \ infinite,}
 
\hfill\textit{but which\ \ just like\ \ the finite\ \ can be determined}
 
\hfill\textit{by definite\ \ well-defined\  and distinguishable num-}
 
\hfill\textit{bers}\ \ [notations].\qquad\qquad\qquad\qquad\qquad\qquad\qquad\qquad
 
\vskip 2 mm

\hfill Georg Cantor, 1883, as translated in \cite{hallett1984} (p. 39)

\vskip 6 mm

Cantor's \lq\lq\textit{unbounded ascending ladder of infinite modes}'' is \emph{sequential} by its very definition -- and, thus, abiding by \emph{the transfinite version} of the Principle of Local Causation, which assumes in this case the form of

\begin{axiom}

\textbf{: Principle of Transfinite Induction.}

Suppose $\Omega$ is a well-ordered totality (set or class), with elements $\alpha,\beta,\ldots$.

$\bullet$ {\bf Induction Hypothesis.} Suppose the formula $\phi(\beta)$ holds for all $\beta<\alpha$.

$\bullet$ {\bf Inductive Step.} Suppose that one can prove that whenever the above \emph{Induction Hypothesis} is verified, the formula $\phi(\alpha)$ holds as well.

$\bullet$ {\bf Inductive Conclusion.} Then $\phi$ holds for all elements of $\Omega$. 

\end{axiom}

\noindent
(Cf., e.g., the discussion of the Induction principle by Joseph R. Shoenfield in \cite{barwise1982} (p. 332))

\subsection{Non-Local Causation of Quantum Algorithms}\label{quantum}

With the discovery of physical, quantum-theoretical non-locality and \emph{quantum entanglement} (from $\bm{EPR}$ paradox,1935, to Bell's inequalities, 1964, to Aspect's experimental confirmation 1982) and its applicability to \emph{QIP}, \emph{non-local causation} of quantum entanglements became the source of dramatic improvements in the efficiency of algorithmic solutions of some important problems, such for example, as the factorization of natural numbers \cite{nielsen_chuang2000}.

The radical novelty of \emph{quantum algorithms} is their \emph{high intrinsic parallelism}, in contradistinction to the \emph{extrinsic, limited, ad hoc parallelism} of classical algorithms. 

In other words, the parallelism of a particular classical \emph{deterministic} algorithm (deploying this suggestive foundational analogy, we do not need to complicate matters by considering {non-deterministic} algorithms) could be achieved only at the second stage of its construction, the first stage being its \emph{strictly sequential} conception and realization of all of its steps -- bearing on, mirroring, and mimicking deterministic processes of classical physics. To this should be added that the functioning of a classical computer executing such an algorithm could be checked, at least in principle, at any moment and in any of its physical bits without interrupting its calculations.

As to the functioning of a quantum computer (as yet, at least materially, a sci-fi dream), it bears on, mirrors, and mimics quantum-theoretical processes, and any attempt to check the state of one of its \emph{qubits} would irreparably corrupt calculations. A quantum algorithm, as its classical companion, is built in \emph{basic sequential stages}, or routines, but, in the difference of a classical algorithm, it does not program the execution of a stage when it represents a genuine quantum process playing out on a carefully chosen bunch of qubits. 

In other words, arriving at such a stage, the algorithm programs the initial tuple of quantum states of the bunch and directs the \lq\lq \emph{quantum flow}'' of information through a carefully chosen sequence of quantum operators -- exactly in the same way as one might direct the flow of water in an irrigation system, without being aware of, interested in, and being able to check what might happen with a particular molecule or bigger portion of water. 

At the end of this loosely controlled flow of parallel processes, the algorithm is supposed to check one or more numeric characteristics of some of the qubits involved (quantum operation of \emph{measuring}), thus destroying the ongoing process, but collecting pieces of information which are used either as a part of the output, or at next quantum stages, or for adjusting the global strategy, or for a combination of all these purposes.

\subsection{Fundamental Object-Subject Duality of Mathematics. Local Causation Case}\label{duality1}

It would be fair to acknowledge that, at the present juncture, the prospects of a radical improvement of the strength, flexibility, universality, and user-friendliness of algorithmic methods of \emph{QIP}, not to speak of the challenges of their hardware realizations, are seriously hampered by the extreme narrowness of the only open today window of genuinely quantum-theoretical opportunities: the entanglement phenomenon. 

The current state of \emph{QIP} affairs should not prevent us from fully appreciate an epistemic and explanatory potential of the non-locality concept -- the potential unimpeded by the currently perceived limits of quantum algorithmic non-locality and bearing on a mysterious duality between non-discrete, non-local continuum phenomena in (logically unrestricted) Mathematics and non-discrete and non-local nature of (proof-theoretically unrestricted) \emph{Mathematical Reasoning}.

In particular, there is no doubt that the Constructivist grievances expressed above so eloquently by the late Errett Bishop (the epigraph to Chapter~\ref{locality}) have something to do with the \emph{exclusively sequential, deterministic, and local} character of mathematical reasoning and verification accepted by the Constructivist school. -- Whereas a \lq\lq\emph{Classical mathematician}'' permits himself to freely \emph{fly over}, or \emph{surf along the continuum flow}, and to do this in a \lq\lq formally irresponsible'', i.e., essentially non-local, non-sequential, and non-deterministic from the Aristotelian logics point of view, way.

We will approach this hypothetic general duality from below, in three \emph{local causation} steps, starting with the most primitive, and yet already very robust and highly important level of Mathematical Reasoning -- \emph{Computation}, the Mathematics 
\lq\lq\textit{that can be carried out by finite beings}'' (Errett Bishop, cf. the epigraph above) named robots, or computers.
\begin{thesis}

\textbf{First Duality Principle: Finitist Mathematics Abides by Local Causation.}

$\bullet$ Constructivism is the branch of Mathematics with the vocation of designing constructive definitions (of creation) of finite mathematical objects and of operations on such objects. Constructive methods are finitist, discrete, and of local causation nature. 

$\bullet$ Classical Theory of Computing (Turing Machine, Church's $\lambda-$calculus, etc.) has the vocation of empowering people with formal finite and discrete computational mechanisms abiding by the condition of local causation, to manipulate objects of constructive nature using constructive operations.

\end{thesis}

The next level of our universal duality thesis corresponds more or less to Errett Bishop's dream of \lq\lq\textit{man's mathematics}'', alias Stephen G. Simpson's aforementioned \lq\lq \textit{countable mathematics}'' (\cite{simpson1999}, p. 1). 

We speak here about problem solving and proof-theoretic analysis which could be carried out in subsystems of Second Order Arithmetic $\bm{Z}_2$ \cite{simpson1999} under the umbrella of the Countable Ordinal Analysis which employs Ordinal Arithmetic on countable ordinals. Ordinal Arithmetic on countable ordinals has the primary purpose of designing systems of ordinal notations to name each ordinal in certain initial stretch of the countable ordinals (cf., e.g., \cite{simmons2004}, p. 65). According to the above analysis, Section~\ref{ordinal_locality}, Countable Ordinal Arithmetic is of a demonstratively local causation nature.

\begin{thesis}

\textbf{Second Duality Principle: Countable Mathematics Abides By Local Causation.}

$\bullet$ The objects of Countable Mathematics are number-theoretical, combinatorial, and analytic theories constructed from axioms of good-behaving subsystems of Second Order Arithmetic, as well propositions, provable in such theories.

$\bullet$ Proof-theoretic machinery employing Countable Ordinal Analysis has the vocation of empowering people with means to both analyze the ordinal rank or complexity of such theories and propositions (properly Ordinal analysis \cite{rathjen2006}) and to determine the minimal in its ordinal complexity set of axioms needed to proof a particular proposition (Reverse Mathematics \cite{simmons2004}).

\end{thesis}

The last \emph{locally causal} level of our universal duality thesis should correspond to the case of {Predicativity}. In one of its completely settled interpretations, {the predicativity attribute} belongs to Countable Mathematics, being of the Feferman-Sch\"utte ordinal rank $\Gamma_0$ (called therefore \emph{the least non-predicatively provable countable ordinal}), and thus, belongs to our Second Duality Thesis. 

However, taking in account the fact that the intuitive notion of predicativity is far from being settled, persisting to remain informal with Solomon Ferferman's last \emph{unfolding} concept, we will restrict ourselves here to the remark that the very attribute of \emph{unfolding} (\cite{feferman2005}, p. 614) perfectly fits in our interpretation of \emph{local causation}.

\section{Non-Locality of Classical Mathematics \emph{Versus} Local Causation of \lq\lq Man's Mathematics''}\label{non-locality}

\subsection{Making the Case For a Synthesis of the Realist and Platonic Conceptions of Mathematics}\label{synthesis}

We are prepared now to treat a conception of Mathematics fitting in Georg Cantor's original, optimistic epistemic vision, the conception both most general from the mathematical practitioner's point of view and most generous from the point of view of a logician. 

The epistemic framework we adopt here is both Platonic whenever we deal with the objectivity and striking relevancy -- or, in Eugene Paul Wigner's words, \textit{\lq\lq the unreasonable effectiveness''} \cite{wigner1960} -- of the flow of the mathematical knowledge in natural sciences, and Realist whenever we need to treat the question of the actual existence of mathematical \lq\lq objects'' whose notations and properties are the real, and only, subjects of any formal mathematical study.

In other words, we do not follow the unabashed, extrovert Platonism of Georg Cantor and, to some degree, of Kurt G\"odel so far as, e.g., to recognize sets as having independent of human definitions existence. And yet, we do assume that, say, the two millenarian quest for the mystery of the real continuum has the merit to bring off correct and pithy presentations of some objectively existing aspects of \lq\lq the real world'', the aspects which, prior to any formalization, are the subject of our mathematical intuition. 

Still, those will be not these aspects which we are attempting to explicitly accommodate in our epistemic and formal framework, but the free flow of mathematically efficacious and formally attractive constructions and theories inspired by this intuition -- and we are intending to carry on this accommodation without any logical or philosophical preconceptions of the last century's metamathematical schools.

We believe that in choosing this middle way between the Scylla of a free wheeling Idealism and Charybdis of a disciplinarian's Pragmatism, we remain faithful to Georg Cantor's original intuition of Mathematical Infinity, all his \emph{post-factum} Platonic justifications notwithstanding. 

As, for example, in the following famous, theologically colored remark, Cantor actually defends not the objective, independent of us existence of transfinite totalities but the free, unimpeded \emph{modus operandi} of the Creator of the Universe -- and thus, \emph{nolens volens}, our own intellectual freedom -- to discover the formal traces of what is captured by Mathematical Infinity and to treat them where, when, and however it might be worth one's while: 

\lq\lq\textit{I am so in favor of the actual infinite that instead of admitting that Nature abhors it, as is commonly said, I hold that Nature makes frequent use of it everywhere, in order to show more effectively the perfections of its Author.}'' (As quoted in \cite{dauben1979}, p. 124.)

\subsection{The Object-Subject Duality of Mathematics. Non-Locality Case: Heuristics}\label{duality}

Contrary to the above duality theses, Section~\ref{duality1}, our choices of \lq\lq Object'' and \lq\lq Subject'' of Mathematics abiding by non-locality causation cannot be as specific, and the corresponding non-locality thesis will remain heuristic. It will be soon clear why it should be so and what it would take to formally advance our understanding of, and mastering Mathematics free from locality constrains.

$\bullet$ \emph{An Instance of a Mathematical Object Abiding by Non-Local Causation.} To begin with, consider the setting of Zeno's paradox \lq\lq Achilles and the Tortoise'' consisting of two independent, simultaneously unfolding, \lq\lq objective'' mathematical processes: the Tortoise's steady advancement along the real line $\mathbb{R}$ and Achilles ever shortened strides along the same line. Achilles strategy and its realization are familiar to us: they are \emph{discrete, recursive, abiding by local causation}. 

But what about the Tortoise's ? Why is Zeno suggesting that hers is a more enigmatic logical and mathematical enterprise, the enterprise incompatible or difficultly compatible with \lq\lq\textit{man's mathematics}'' of Errett Bishop ? -- By the way, it is not by chance that Zeno has chosen a \emph{man}, Achilles, who has to catch up with a fabulously slow and yet mysteriously unattainable Tortoise. 

The answer, we believe, is that Zeno allots to the Tortoise the \lq\lq\textit{divine ability}'' to \lq\lq calmly swim along the flow'' of the continuum, whereas Achilles, the man, has no choice but to dangerously hang over the bottomless depths of this primeval continuum stream, being able to advance only by jumping from a one point-size location on this stream -- the location carefully chosen and temporally frozen solid -- to the next one. 

In other words, Zeno's paradox suggests, first, that the real continuum $\mathbb{R}$ -- at any given moment of our considerations -- \lq\lq divinely'' (taking Errett Bishop's at his word) \emph{exists in our intuition}, simultaneously with all its continuum points, intervals, etc. -- the scenario clearly implying some \emph{non-local causation} -- and second, that any specific mathematical question implying $\mathbb{R}$ cannot be \emph{humanly resolved} but within an appropriate discrete, sequential mathematical framework abiding by \emph{a local causation principle}.

So much about the aspects of non-local causation in our \emph{intuitive perception} of the real continuum. 

\subsection{Non-Locality. Interpretation \& Foundational Implications I: Logical Paradoxes, Circularity, and Truth}\label{circularity}

Moreover, we claim that the discovery, by Bertrand Russel and his followers, of logical paradoxes characterized by \emph{self-referentiality}, or as it is known today \emph{circularity} of formal arguments \cite{finsler1996} -- such as Russel's \lq\lq\textit{set of all sets which are not elements of themselves}'' or Liar paradoxes \cite{martin1984} -- were, and to some extent, remain just new, surprising, and to the pioneers and early practitioners of our science, historically painful manifestations of the authenticity and fundamentality of Zeno's and Everrett Bishop's predicament: how to represent an essentially and irremediably non-local \lq\lq divine'' reality by \lq\lq man's'' sequential local causality arguments.

\begin{making}

\textbf{: Two Types of Logical Paradoxes Reflecting Two Sources of Uncertainty about the Truth of Mathematical and Logical Claims.}

{\bf (1)} Russel's paradox simulates the unavoidable circularity of some definitions of set theory and of the original, full-blown, non-constructive and non-predicative system of analysis. As Zeno's paradox, it pushes to its absurd limits a conventional wisdom about the human ability to \emph{instantly} and yet correctly appreciate a complex, dynamical, nonlocal by its very nature reality. Whereas in Zeno's case Achilles possesses the faculty of an \emph{instant observation} of the Tortoise's position but is artificially deprived of the faculty of \emph{instant prevision} concerning her future positions, Russel's paradox plays in a similar way with human faculties to be \emph{instantly formally correct} in a very selectively, syntactically understood way -- but to pitifully fail semantically.

{\bf (2)} As to Liar paradoxes, from simple to strengthened ones, they play with the natural languages ability to treat nonlocal reality, nonlocal both \emph{temporally} -- i.e., outside the limits of \emph{instantaneousness} -- and spatially but also relationally, in a formally ambiguous but perfectly understandable and workable way. 

\end{making}
 
\subsection{Non-Locality. Interpretation \& Foundational Implications II: Problem Solving, Theorem Proving, Theories Building}\label{building}

$\bullet$ \emph{Mathematical Subject Abiding by Non-Local Causation.} As a matter of course, similar aspects of non-locality are present in the treatment of all cases of infinite, or even finite but immensely big mathematical structures. In fact a good practitioner should always have a global, all-embracing, instantaneous, non-local in its causation, intuitive grasp of structures of interest to him, both exemplified and pushed to its limits by Srinivasa Ramanujan's (1887-1920) \emph{intimate familiarity with numbers} (\cite{hardy1940}, p. 12) and his general way of doing mathematics -- by \lq\lq\textit{showing astonishing imaginative power}'', even if \lq\lq\textit{as always, proving next to nothing}'' \cite{hardy1940} (p. 15). Two other young and tragically early departed mathematical geniuses, Niels Henrik Abel (1802-1829) and Evarist Galois (1811-1832), have left to the posterity not less rich and original visions of both never before observed mathematical landscapes and mathematical laws by which these landscapes abide -- but not much, or not at all, traces of deductive rigor and very little, if at all, proofs.

From the epistemic level of \emph{Computation}, we have already ascended -- up the ladder of epistemic non-triviality and intrinsic conceptual complexity of mathematical challenges -- to \emph{Theorem Proving}, \emph{Axiomatizing}, \emph{Consistency, Independency, and Ordinal Strength} of axiomatic studies. Doing this, we were following the metamathematical practices of the 20th century which have successfully formalized, customized, and perfected computational and proof-theoretic methods and skills -- always abiding, as it has been demonstrated above, Chapter~\ref{locality}, by local causation principles.

The next step should dispense us with this fundamental restriction and to bring first, still tentative, fuzzy, and yet indubitable insights into what might one day become \emph{Mathematics with Non-Local Causation}. 

There is not much precedents for, and apparently no precursors of such an enterprise. The following historical remark by Igor R. Shafarevich gives the closest known to the present author fit for an extrinsic -- as antonymous to \emph{intrinsic}, in other words, epistemic, i.e., of commanding interest to the present study -- foretaste of the \emph{enigma of non-local causation}: 

\lq\lq\textit{Viewed superficially, Mathematics is the result of centuries of effort by many thousands of largely unconnected individuals scattered across continents, centuries and millennia. However, the internal logic of its development much more resembles the work of a single intellect developing its thought in a continuous and systematic way, and only using as a means a multiplicity of human individualities, much as in an orchestra playing a symphony written by some composer the theme moves from one instrument to another so that as soon as one performer is forced to cut short his part, it is taken up by another player, who continues it with due attention to the score.}'' \cite{shafarevich1981} (p. 182) 

Without pretending to fully understand, let alone to endorse the presumed by the above statement \emph{providential character} of the progress in Mathematics driven by apparently unrelated efforts of a dispersed and disconnected collectivity of professionals, we choose to apply here Shafarevich's noble vision of \emph{new, emerging mathematical artifacts and theories} to the work of an individual: as a simultaneous, harmoniously coordinated execution of, say, a symphony by \lq\lq a band of players'', i.e., her or his disparate, almost unrelated, in many cases most illogical insights \lq\lq playing the whole thing out'' without any restraints of local causation.

Moreover, we see the authentic vocation of \emph{Non-Local Causation Mathematics} as 
 \emph{Metamathematics of Emerging Mathematical Objects, Claims, Theories}:
\begin{thesis}

\textbf{Third Duality Principle: Non-Local Causation in Problem Solving, Theorem Proving, Theories Building.}

With the emergence of non-local causation in the treatment of the aforementioned challenges, we are formally approaching the next level of mathematical reasoning: the emergence, sometimes, simultaneously or in a casual order, first, of initially only intuitively perceived new mathematical structures, accompanied by an intuitive conception of eventual claims concerning these structures, and followed by an intuitive search for proofs of such claims. 

\end{thesis}

{}

\noindent
\textit{Edward G. BELAGA}

\noindent
\textit{Institut de Recherche Math\'{e}matique Avanc\'{e}e de Strasbourg} 

\noindent
\textit{7 rue
Ren\'{e} Descartes, F-67084 Strasbourg Cedex, FRANCE}

\noindent
e-mail: \textit{edward.belaga@math.u-strasbg.fr}

\end{document}